\documentclass[11pt]{article}
\usepackage{geometry}                
\usepackage{amsmath,amssymb}
\usepackage{amsthm}
\usepackage{hyperref}
\usepackage{graphicx,epsfig,color}
\usepackage{booktabs,bm,multirow}
\usepackage{cases}
\usepackage[caption=false]{subfig}
\newtheorem{theorem}{Theorem}
\newtheorem{proposition}[theorem]{Proposition}
\newtheorem{remark}[theorem]{Remark}
\newtheorem{lemma}[theorem]{Lemma}

\def\mbbr{\mathbb R}
\def\mcalk{\mathcal K}
 
\def\l{\left} 
\def\r{\right} 
\def\mbf{\mathbf}

\def\rmi{{\rm i}}

\def\spa{{\rm span}}
\def\rank{{\rm rank}} 
\def\ran{{\rm range}} 
\def\nul{{\rm null}}
\def\diag{{\rm diag}}

\DeclareMathOperator*{\argmin}{argmin}

\def\dsp{\displaystyle} 

\def\wt{\widetilde}

\newcommand{\beq}{\begin{equation}}\newcommand{\eeq}{\end{equation}}
\newcommand{\bit}{\begin{itemize}}\newcommand{\eit}{\end{itemize}}
\newcommand{\bem}{\begin{bmatrix}}\newcommand{\eem}{\end{bmatrix}}

\title{On Krylov subspace methods for skew-symmetric and shifted skew-symmetric linear systems}

\author{Kui Du\thanks{School of Mathematical Sciences and Fujian Provincial Key Laboratory of Mathematical Modeling and High Performance Scientific Computing, Xiamen University, Xiamen 361005, China (kuidu@xmu.edu.cn).},\quad  Jia-Jun Fan\thanks{School of Mathematical Sciences, Xiamen University, Xiamen 361005, China ({fanjiajun@stu.xmu.edu.cn}).},\quad  Xiao-Hui Sun\thanks{School of Mathematics and Statistics, Xiamen University of Technology, Xiamen 361024, China ({2023000058@xmut.edu.cn}).},\quad  Fang Wang\thanks{School of Mathematical Sciences, Xiamen University, Xiamen 361005, China ({fangwang@stu.xmu.edu.cn}).},\quad  Ya-Lan Zhang\thanks{School of Mathematical Sciences, Xiamen University, Xiamen 361005, China ({zhangyalan@stu.xmu.edu.cn}).}} 
\date{}                                           

\begin{document}
\maketitle

\begin{abstract}  
\vspace{.5mm} 

Krylov subspace methods for solving linear systems of equations involving skew-symmetric matrices have gained recent attention. Numerical equivalences among Krylov subspace methods for nonsingular skew-symmetric linear systems have been given in Greif et al. [SIAM J. Matrix Anal. Appl., 37 (2016), pp. 1071--1087]. In this work, we extend the results of Greif et al. to singular skew-symmetric linear systems. In addition, we systematically study three Krylov subspace methods (called S$^3$CG, S$^3$MR, and S$^3$LQ) for solving shifted skew-symmetric linear systems. They all are based on Lanczos triangularization for skew-symmetric matrices, and correspond to CG, MINRES, and SYMMLQ for solving symmetric linear systems, respectively. To the best of our knowledge, this is the first work that studies S$^3$LQ. We give some new theoretical results on S$^3$CG, S$^3$MR, and S$^3$LQ. We also provide the relationship among the three methods and those based on Golub--Kahan bidiagonalization and Saunders--Simon--Yip tridiagonalization. Numerical examples are given to illustrate our theoretical findings.



\end{abstract}

\section{Introduction} 
We are concerned with Krylov subspaces methods for solving a linear system $$\mbf A\mbf x=\mbf b,$$
where $\mbf A\in\mbbr^{n\times n}$ is a skew-symmetric (i.e., $\mbf A^\top=-\mbf A$) or shifted skew-symmetric (i.e., $\mbf A=\alpha\mbf I+\mbf S$, $\alpha\neq 0$, $\mbf S^\top=-\mbf S$) matrix, $\mbf b\in\mbbr^n$, $\mbf A\neq \mbf 0$, and $\mbf b\neq \mbf 0$. Such systems arise in numerous applications, including computational fluid dynamics, optimization, and systems theory, etc.; see \cite{bai2003hermi,idema2007minim,greif2009itera,chen2023accel,idema2023compa} and the references therein. Let $\mbf A^\dag$ denote the Moore--Penrose pseudoinverse of $\mbf A$.  When $\mbf A$ is singular and skew-symmetric, we are interested in the {\it minimum-length} or {\it pseudoinverse} solution $\mbf A^\dag\mbf b$, which is the unique solution of the related least-squares problem $$\min\|\mbf x\|_2\quad  \text{ subject to }\quad \mbf x\in\argmin_{\mbf x\in\mbbr^n}\|\mbf b-\mbf A\mbf x\|_2.$$

 Greif and Varah \cite{greif2009itera} proposed two Krylov subspace methods for solving {\it nonsingular} skew-symmetric linear systems, one based on the Galerkin condition, the other on the minimum residual condition. For the sake of convenience, we will refer to these two methods as S$^2$CG and S$^2$MR. Eisenstat \cite{eisenstat2015equiv} gave algorithm-independent proofs that in exact arithmetic the iterates for S$^2$CG and S$^2$MR are identical to the iterates for the classic Craig's method \cite{craig1955nstep} (also called CGNE \cite[Chapter 7]{greenbaum1997itera}) and CGLS \cite{paige1982lsqr} (also called CGNR \cite[Chapter 7]{greenbaum1997itera}), respectively. By using the numerical equivalence of Lanczos tridiagonalization  \cite{lanczos1950itera} and Golub--Kahan bidiagonalization \cite{golub1965calcu} for skew-symmetric $\mbf A$ and nonzero $\mbf b$, Greif et al. \cite{greif2016numer} showed that S$^2$CG and S$^2$MR are numerically equivalent to Golub--Kahan bidiagonalization variants of CGNE and CGNR, respectively. We note that there is no literature on the performance of S$^2$CG and S$^2$MR applied to {\it singular} skew-symmetric linear systems. 

Galerkin and minimum residual methods for solving linear system $\bf Ax=b$ with positive definite $\mbf A$ (i.e., the symmetric part $(\mbf A+\mbf A^\top)/2$ is positive definite) have been discussed; see, e.g.,  \cite{concus1976gener,widlund1978lancz,rapoport1978nonli,szyld1993varia,guducu2022nonhe}. It is straightforward to show that these methods can be used to solve shifted skew-symmetric linear systems. Similar to MINRES \cite{paige1975solut} for solving symmetric linear systems, Idema and Vuik \cite{idema2007minim,idema2023compa} and Jiang \cite{jiang2007algor} have introduced a minimum residual method for shifted skew-symmetric linear systems by using Lanczos tridiagonalization. For the sake of convenience, we will refer to their method as S$^3$MR. Similar to the derivation of CG  \cite{hestenes1952metho} for symmetric positive definite linear systems (see \cite[Chapter 6]{demmel1997appli}), one can design a Galerkin method for shifted skew-symmetric linear systems. We call the resulting method S$^3$CG. As we all know, SYMMLQ \cite{paige1975solut} for solving symmetric linear systems can produce an approximate solution sequence with decreasing error norms; see \cite{szyld1993varia}. At present, there does not exist a method corresponding to SYMMLQ to solve shifted skew-symmetric linear systems in the literature. Here we develop such a solver and call it S$^3$LQ.   

The main contributions of this work are as follows. (i) We extend the results of Greif et al. \cite{greif2016numer} to singular skew-symmetric linear systems. Specifically, we show that S$^2$CG returns the pseudoinverse solution for singular consistent skew-symmetric linear systems, and S$^2$MR returns the pseudoinverse solution for arbitrary singular skew-symmetric linear systems (including consistent and inconsistent cases). (ii) We develop a new solver, S$^3$LQ, which completes the family of Krylov subspace methods for shifted skew-symmetric linear systems based on Lanczos tridiagonalization and is a natural companion to S$^3$CG and S$^3$MR. (iii) We provide the relation among S$^3$LQ, S$^3$CG, and Craig's method, and the relation between S$^3$MR and LSQR. (iv) We point out the relation between Lanczos tridiagonalization and Saunders--Simon--Yip tridiagonalization \cite{saunders1988two} for shifted skew-symmetric matrices. Based on this relation, we show that S$^3$LQ and S$^3$MR coincide with USYMLQ and USYMQR \cite{saunders1988two}, respectively, when applied to shifted skew-symmetric linear systems.

The paper is organized as follows. In the rest of this section, we give other related research and some notation. Section 2 gives background on Lanczos tridiagonalization, Golub--Kahan bidiagonalization, and Saunders--Simon--Yip tridiagonalization. Sections 3 and 4 discuss Krylov subspace methods for skew-symmetric and shifted skew-symmetric linear systems, respectively. Numerical experiments are given in section 5, followed by concluding remarks in section 6.

{\it Related research}. For skew-symmetric linear systems, Huang et al. \cite{huang1999itera} gave an iterative method computing only the even minimum residual iterates. Idema and Vuik \cite{idema2007minim} showed that the method in \cite{huang1999itera} is equivalent to CGNR. Gu and Qian \cite{gu2009skew} imposed the Galerkin condition and presented skew-symmetric methods for solving nonsymmetric linear systems with multiple right-hand sides. Chen and Shen \cite{chen2007regul} proposed a regularized conjugate gradient method for solving ill-conditioned skew-symmetric linear systems. They also showed that their method can be used to solve ill-conditioned shifted skew-symmetric or ill-conditioned skew-symmetrizable linear systems. Recently, Chen and Wei \cite[section 3.4]{chen2023accel} proposed accelerated overrelaxation methods for shifted skew-symmetric linear systems, and showed their linear convergence by using a Lyapunov analysis.

There are some studies on Krylov subspace methods for singular linear systems in the literature; see, for example, \cite{brown1997gmres,ipsen1998idea,smoch1999some,sidi1999unifi,wei2000conve,calvetti2000gmres,cao2002note,reichel2005break,smoch2007spect,hayami2011geome,choi2011minre,choi2013minim,morikuni2018gmres}. For singular $\mbf A$, conditions under which the GMRES \cite{saad1986gmres} iterates converge safely to a least-squares solution or to the pseudoinverse solution are given in \cite{brown1997gmres}. In particular, the pseudoinverse solution is guaranteed (see \cite[Theorem 2.4]{brown1997gmres}) when $\mbf A$ is range-symmetric, the initial approximation is in the range of $\mbf A$, and the system is consistent. Choi \cite{choi2013minim} proposed SS-MINRES-QLP for solving singular skew-symmetric linear systems, and showed that whether the systems are consistent or not, SS-MINRES-QLP computes the pseudoinverse solution. But she did not consider S$^2$MR.

{\it Notation}. Lowercase (uppercase) boldface letters are reserved for column vectors (matrices). Let $\rmi=\sqrt{-1}$ be the imaginary unit. For any vector $\mbf b \in\mbbr^n$, we use $\mbf b^\top$ and $\|\mbf b\|_2$ to denote the transpose and the Euclidean norm of $\mbf b$, respectively. We use $\mbf I_k$ to denote the $k\times k$ identity matrix, and use $\mbf e_k$ to denote the $k$th column of the identity matrix whose order is clear from the context. For any matrix $\mbf A \in\mbbr^{m\times n}$, we use $\mbf A^\top$, $\mbf A^\dag$, $\ran(\mbf A)$, $\nul(\mbf A)$, and $\rank(\mbf A)$ to denote the transpose, the Moore--Penrose pseudoinverse, the column space, the null space, and the rank of $\mbf A$, respectively. For nonsingular $\mbf A$, we use $\mbf A^{-1}$ to denote its inverse. Any normalization of the form $\mbf \gamma \mbf w:=\mbf b$ is short for ``$\gamma:=\|\mbf  b\|_2$; if $\gamma=0$ then stop; else $\mbf w=\mbf b/\gamma$.''

\section{Preliminaries}

\subsection{Lanczos tridiagonalization for skew-symmetric matrices.}

Given a nonzero vector $\mbf b\in\mbbr^n$, the Arnoldi iteration \cite{arnoldi1951princ} for a general matrix $\mbf  A\in\mbbr^{n\times n}$ generates a sequence of orthonormal vectors $\{\mbf w_k\}$ such that $\mbf w_1=\mbf b/\gamma_1$ with $\gamma_1=\|\mbf b\|_2$, $\mbf W_k^\top\mbf W_k=\mbf I_k$, and $$\mbf A\mbf W_k=\mbf W_{k+1}\mbf H_{k+1,k}=\mbf W_k\mbf H_k+\gamma_{k+1}\mbf w_{k+1}\mbf e_k^\top, \quad \mbf H_k=\mbf W_k^\top\mbf A\mbf W_k,$$ where $\mbf W_k:=\bem \mbf w_1 &\mbf w_2 & \cdots & \mbf w_k \eem$, $\mbf H_{k+1,k}$ is a $(k+1)\times k$ upper Hessenberg matrix, and $$\ran(\mbf W_k)=\mcalk_k(\mbf A,\mbf b):=\spa\{\mbf b, \mbf A \mbf b, \cdots, \mbf A^{k-1}\mbf b\}.$$ If $\mbf A^\top=-\mbf A$, then $\mbf H_k^\top=\mbf W_k^\top\mbf A^\top\mbf W_k=-\mbf W_k^\top\mbf A\mbf W_k=-\mbf H_k$, giving $$\mbf H_{k+1,k}=\bem 0 & -\gamma_2 &&\\ \gamma_2 & 0 & \ddots& \\ &\ddots &\ddots &-\gamma_k \\ &&\gamma_k & 0 \\ &&&\gamma_{k+1} \eem=\bem \mbf H_k\\ \gamma_{k+1}\mbf e_k^\top \eem.$$ This leads to Algorithm 1, a short recurrence Lanczos process computing $\{\mbf w_k\}$ and $\{\gamma_k\}$. Note that in Algorithm 1 each $\gamma_k$ results simply from the normalization of its $\mbf w_k$.
\begin{center}
\begin{tabular*}{170mm}{l}
\toprule {\bf Algorithm 1}: Lanczos tridiagonalization for skew-symmetric $\mbf A\in\mbbr^{n\times n}$ and nonzero $\mbf b\in\mbbr^n$ 
\\ \hline\noalign{\smallskip}
\qquad Compute $\gamma_1\mbf w_1:=\mbf b$ and $\gamma_2\mbf w_2:=\mbf A\mbf w_1$. \\\noalign{\smallskip}
\qquad {\bf for} $k=2,3,\ldots$ {\bf do}\\ \noalign{\smallskip}
\qquad\qquad $\gamma_{k+1}\mbf w_{k+1}:=\mbf A\mbf w_k+\gamma_k\mbf  w_{k-1}$;\\ \noalign{\smallskip}
\qquad {\bf end}\\ 
\bottomrule
\end{tabular*}
\end{center}

There are at most $n$ orthonormal vectors, so theoretically Algorithm 1 must stop in $\ell\leq n$ steps. The positive integer $\ell$ satisfies $$\dim\mcalk_k(\mbf A,\mbf b)=\begin{cases}
k, & \text{ if } 1\leq k\leq \ell \\
\ell, & \text{ if } k\geq\ell+1. 
\end{cases}$$  At the $\ell$th step, we have $\gamma_k>0$ for $k=1,2\ldots,\ell$, $\gamma_{\ell+1}=0$, $\mbf w_{\ell+1}=\mbf 0$, and $$\mbf W_\ell^\top\mbf W_\ell=\mbf I_\ell,\quad\mbf A\mbf W_\ell=\mbf W_\ell\mbf H_\ell,\qquad \mbf H_\ell=\mbf W_\ell^\top\mbf A\mbf W_\ell=\bem 0 & -\gamma_2 &&\\ \gamma_2 & 0 & \ddots &  \\ &\ddots &\ddots &-\gamma_\ell  \\ &&\gamma_\ell & 0 \eem.$$ We provide some properties of $\mbf H_\ell$ and its submatrices in the following theorem. We note that a similar result has been presented in \cite[Theorem 2.1]{choi2013minim}.

\begin{theorem}\label{oddeven} Assume that $\mbf A^\top=-\mbf A$. For each $j$ with $1\leq j\leq\ell/2$,  $\mbf H_{2j}$ is nonsingular. If $\mbf b\in\ran(\mbf A)$, then $\ell$ is even and $\mbf H_\ell$ is nonsingular. Otherwise, $\ell$ is odd and $\mbf H_\ell$ is singular. 
\end{theorem}
\proof  
Since $\mbf H_k$ is real skew-symmetric, it is unitarily diagonalizable, and its non-zero eigenvalues are pure imaginary complex conjugate pairs. This means that $\rank(\mbf H_k)$ is even. For each $k$ with $2\leq k\leq \ell$, the $(k-1)\times(k-1)$ lower-left submatrix of $\mbf H_k$ is nonsingular, which implies $\rank(\mbf H_k)\geq k-1$. Therefore, for each $j$ with $1\leq j\leq\ell/2$, we have $\rank(\mbf H_{2j})=2j$, i.e., $\mbf H_{2j}$ is nonsingular.

If $\mbf b\in\ran(\mbf A)$, then $\mbf w_1\in\ran(\mbf A)$. Now, assume that $\ell$ is odd. From $\gamma_{k+1}\mbf w_{k+1}=\mbf A\mbf w_k+\gamma_k\mbf  w_{k-1}$, we have $\mbf w_{2j-1}\in\ran(\mbf A)$ for each $j$ with $1\leq j\leq (\ell+1)/2$. Since $\mbf H_\ell$ is singular, there exists a nonzero vector $\mbf y\in\mbbr^\ell$ satisfying $\mbf H_\ell\mbf y=\mbf 0$. It is straightforward to show that even entries of $\mbf y$ must be zero, i.e., $\mbf y=\bem y_1 & 0 & y_3 & 0 & \cdots  & y_\ell\eem^\top$. Hence, $\mbf W_\ell\mbf y\in\ran(\mbf A)$. By $\mbf A^\top=-\mbf A$, $\mbf A\mbf W_\ell=\mbf W_\ell\mbf H_\ell$, and $\mbf H_\ell\mbf y=\mbf 0$, we have $\mbf A^\top\mbf W_\ell\mbf y=-\mbf A\mbf W_\ell\mbf y=-\mbf W_\ell\mbf H_\ell\mbf y=\mbf 0$,  which gives $\mbf W_\ell\mbf y\in\nul(\mbf A^\top)$. By $\ran(\mbf A)\cap\nul(\mbf A^\top)=\{\mbf 0\}$, we have $\mbf W_\ell\mbf y=\mbf 0$, which yields $\mbf y=\mbf 0$. This, however, is a contradiction to $\mbf y\neq\mbf 0$. Therefore, $\ell$ must be even.

If $\mbf b\notin\ran(\mbf A)$, then $\mbf w_1\notin\ran(\mbf A)$. Now, assume that $\ell$ is even. It follows that $\mbf H_\ell$ is nonsingular. By $\mbf A\mbf W_\ell=\mbf W_\ell\mbf H_\ell$, we have $\mbf W_\ell=\mbf A\mbf W_\ell\mbf H_\ell^{-1}$. This means $\mbf w_k\in\ran(\mbf A)$ for each $k$ with $1\leq k\leq\ell$. This, however, is a contradiction to $\mbf w_1\notin\ran(\mbf A)$. Therefore, $\ell$ must be odd.
\endproof

\subsection{Golub--Kahan bidiagonalization} 

For a general matrix $\mbf A\in\mbbr^{n\times m}$ and a nonzero vector $\mbf b\in\mbbr^n$, we describe  Golub--Kahan bidiagonalization \cite{golub1965calcu} as Algorithm 2. 

\begin{center}
\begin{tabular*}{170mm}{l}
\toprule {\bf Algorithm 2}: Golub--Kahan bidiagonalization for general $\mbf A\in\mbbr^{n\times m}$ and nonzero $\mbf b\in\mbbr^n$  
\\ \hline\noalign{\smallskip}
\qquad Compute $\beta_1\mbf u_1:=\mbf b$ and $\alpha_1\mbf v_1:=\mbf A^\top\mbf u_1$. \\\noalign{\smallskip}
\qquad {\bf for} $j=1,2,\cdots$  {\bf do}\\ \noalign{\smallskip}
\qquad\qquad $\beta_{j+1}\mbf u_{j+1}:=\mbf A\mbf v_j-\alpha_j\mbf u_j$;\\ \noalign{\smallskip}
\qquad\qquad $\alpha_{j+1}\mbf v_{j+1}:=\mbf A^\top\mbf u_{j+1}-\beta_{j+1}\mbf v_j$;\\ \noalign{\smallskip}
\qquad {\bf end}\\ 
\bottomrule
\end{tabular*}
\end{center}

We define $\mbf U_j:=\bem \mbf u_1 & \mbf u_2 &\cdots & \mbf u_j\eem$, $\mbf V_j:=\bem \mbf v_1 & \mbf v_2 &\cdots & \mbf v_j\eem$, and $$\mbf B_j:=\bem \alpha_1 &&&\\ \beta_2&\alpha_2 & &\\ &\ddots &\ddots &\\ && \beta_j &\alpha_j \eem,\qquad \mbf B_{j+1,j}:=\bem \alpha_1 &&&\\ \beta_2&\alpha_2 &&\\ &\ddots &\ddots &\\ && \beta_j &\alpha_j \\ &&& \beta_{j+1}\eem=\bem\mbf B_j\\ \beta_{j+1}\mbf e_j^\top\eem.$$ After $j$ iterations of Algorithm 2, the following hold: $$\mbf A\mbf V_j=\mbf U_{j+1}\mbf B_{j+1,j}=\mbf U_j\mbf B_j+\beta_{j+1}\mbf u_{j+1}\mbf e_j^\top, \quad \mbf A^\top\mbf U_{j+1}=\mbf V_{j+1}\mbf B_{j+1}^\top=\mbf V_j\mbf B_{j+1,j}^\top+\alpha_{j+1}\mbf v_{j+1}\mbf e_{j+1}^\top,$$ $$\mbf U_j^\top\mbf U_j=\mbf V_j^\top\mbf V_j=\mbf I_j,\quad \ran(\mbf U_j)=\mcalk_j(\mbf A\mbf A^\top,\mbf b),\quad \ran(\mbf V_j)=\mcalk_j(\mbf A^\top\mbf A,\mbf A^\top\mbf b).$$ In exact arithmetic, Algorithm 2 must stop in $k_0\leq\min(m,n)$ steps and either $\alpha_{k_0+1}=0$ or $\beta_{k_0+1}=0$. Paige \cite{paige1974bidia} showed that $\beta_{k_0+1}=0$ if $\mbf b\in\ran(\mbf A)$ and $\alpha_{k_0+1}=0$ if $\mbf b\notin\ran(\mbf A)$.  

Greif et al. \cite{greif2016numer} showed that one step of Algorithm 2 applied to a skew-symmetric matrix $\mbf A$ is equivalent to two steps of Algorithm 1. More precisely, we have the following relations: $$\beta_j=\gamma_{2j-1},\quad \alpha_j=\gamma_{2j},\quad\mbf u_j=(-1)^{j-1}\mbf w_{2j-1},\quad \mbf v_j=(-1)^j\mbf w_{2j},\quad 1\leq j\leq\ell_0:=\lceil {\ell}/{2}\rceil.$$ If $\ell$ is even, then Algorithm 2 terminates with  $j=\ell_0$ and $\beta_{\ell_0+1}=0$. We have $$\mbf A\mbf V_{\ell_0}=\mbf U_{\ell_0}\mbf B_{\ell_0},\qquad \mbf A^\top\mbf U_{\ell_0}=\mbf V_{\ell_0}\mbf B_{\ell_0}^\top.$$ If $\ell$ is odd, then Algorithm 2 terminates with  $j=\ell_0-1$ and $\alpha_{\ell_0}=0$. We have $$\mbf A\mbf V_{\ell_0-1}=\mbf U_{\ell_0}\mbf B_{\ell_0,\ell_0-1},\qquad \mbf A^\top\mbf U_{\ell_0}=\mbf V_{\ell_0-1}\mbf B_{\ell_0,\ell_0-1}^\top.$$

\subsection{Saunders--Simon--Yip tridiagonalization} 

For a general matrix $\mbf A\in\mbbr^{n\times m}$ and two nonzero vectors $\mbf b\in\mbbr^n$ and $\mbf c\in\mbbr^m$, we describe Saunders--Simon--Yip tridiagonalization \cite{saunders1988two} as Algorithm 3.

\begin{center}
\begin{tabular*}{170mm}{l}
\toprule {\bf Algorithm 3}: Saunders--Simon--Yip tridiagonalization for $\mbf A\in\mbbr^{n\times m}$, $\mbf b\in\mbbr^n$, and $\mbf c\in\mbbr^m$
\\ \hline\noalign{\smallskip}
\qquad Set $\wt{\mbf u}_0=\mbf 0$ and $\wt{\mbf v}_0=\mbf 0$. \\
\qquad Compute $\beta_1\wt{\mbf u}_1:=\mbf b$ and $\alpha_1\wt{\mbf v}_1:=\mbf c$. \\\noalign{\smallskip}
\qquad {\bf for} $k=1,2,\cdots$  {\bf do}\\ \noalign{\smallskip}
\qquad\qquad $\mbf q:=\mbf A\wt{\mbf v}_k-\alpha_k\wt{\mbf u}_{k-1}$;\\ \noalign{\smallskip}
\qquad\qquad $\theta_k:=\wt{\mbf u}_k^\top\mbf q$;\\ \noalign{\smallskip}
\qquad\qquad $\beta_{k+1}\wt{\mbf u}_{k+1}:=\mbf q-\theta_k\wt{\mbf u}_k$;\\ \noalign{\smallskip}
\qquad\qquad $\alpha_{k+1}\wt{\mbf v}_{k+1}:=\mbf A^\top\wt{\mbf u}_k-\beta_k\wt{\mbf v}_{k-1}-\theta_k\wt{\mbf v}_k$;\\ \noalign{\smallskip}
\qquad {\bf end}\\ 
\bottomrule
\end{tabular*}
\end{center}

We define $\wt{\mbf U}_k:=\bem\wt{\mbf u}_1 &\wt{\mbf u}_2 &\cdots & \wt{\mbf u}_k\eem$, $\wt{\mbf V}_k:=\bem\wt{\mbf v}_1 &\wt{\mbf v}_2 &\cdots & \wt{\mbf v}_k\eem$, and $$\wt{\mbf H}_k:=\bem \theta_1 & \alpha_2 && \\ \beta_2 & \theta_2 &\ddots & \\ &\ddots &\ddots & \alpha_k\\ &&\beta_k &\theta_k \eem,\quad \wt{\mbf H}_{k+1,k}:=\bem\wt{\mbf H}_k\\ \beta_{k+1}\mbf e_k^\top\eem,\quad \wt{\mbf H}_{k,k+1}:=\bem\wt{\mbf H}_k &\alpha_{k+1}\mbf e_k\eem.$$ After $k$ iterations of Algorithm 3, the following hold: $$\mbf A\wt{\mbf V}_k=\wt{\mbf U}_{k+1}\wt{\mbf H}_{k+1,k}=\wt{\mbf U}_k\wt{\mbf H}_k+\beta_{k+1}\wt{\mbf u}_{k+1}\mbf e_k^\top,$$ $$\mbf A^\top\wt{\mbf U}_k=\wt{\mbf V}_{k+1}\wt{\mbf H}_{k,k+1}^\top=\wt{\mbf V}_k\wt{\mbf H}_k^\top+\alpha_{k+1}\wt{\mbf v}_{k+1}\mbf e_k^\top,$$ $$\wt{\mbf U}_k^\top\wt{\mbf U}_k=\wt{\mbf V}_k^\top\wt{\mbf V}_k = \mbf I_k,\quad \wt{\mbf H}_k=\wt{\mbf U}_k^\top\mbf A\wt{\mbf V}_k.$$ By induction on $k$, it can be shown that Algorithm 3 for the case $\mbf A^\top=-\mbf A$ and $\mbf c=\mbf b$ reduces to Algorithm 1. More precisely, we have the following relations: $$\theta_k=0,\quad \alpha_k=\beta_k=\gamma_k,\quad \wt{\mbf u}_{2j-1}=\wt{\mbf v}_{2j-1}=(-1)^{j+1}\mbf w_{2j-1},\quad \wt{\mbf u}_{2j}=-\wt{\mbf v}_{2j}=(-1)^{j+1}\mbf w_{2j}.$$ After $\ell$ iterations, we have $$\wt{\mbf U}_\ell=\mbf W_\ell\mbf D_\ell,\quad  \wt{\mbf V}_\ell=\mbf W_\ell\wt{\mbf D}_\ell,\quad \wt{\mbf U}_\ell=\wt{\mbf V}_\ell\wt{\mbf D}_\ell\mbf D_\ell,$$ $$\mbf A\wt{\mbf V}_\ell=\wt{\mbf U}_\ell\wt{\mbf H}_\ell,\quad \mbf A^\top\wt{\mbf U}_\ell=\wt{\mbf V}_\ell\wt{\mbf H}_\ell^\top,\quad \wt{\mbf H}_\ell=\mbf D_\ell\mbf H_\ell\wt{\mbf D}_\ell,$$ where $$\mbf D_\ell=\diag\{1,1,-1,-1,1,1,-1,-1,\cdots,(-1)^{\lfloor\frac{\ell+3}{2}\rfloor}\},$$ and $$\wt{\mbf D}_\ell=\diag\{1,-1,-1,1,1,-1,-1,1,\cdots,(-1)^{\lceil\frac{\ell+3}{2}\rceil}\}.$$

\section{Krylov subspace methods for skew-symmetric linear systems}
In this section, we consider Krylov subspaces methods for skew-symmetric linear systems. We note that the existing results \cite{greif2009itera,eisenstat2015equiv,greif2016numer} are mainly on nonsingular cases. Here we extend these results to singular cases. Throughout the paper, we use the zero vector as the initial approximation for simplicity.  

\subsection{S$^2$CG and CRAIG} 

After $k$ iterations of Lanczos tridiagonalization (see section 2.1) for skew-symmetric $\mbf A$ and nonzero $\mbf b$, the Galerkin condition for solving $\bf Ax=b$ is $$\mbf x_k^{\rm G}=\mbf W_k\mbf y_k^{\rm G},\qquad \mbf b-\mbf A\mbf x_k^{\rm G} \perp \mcalk_k(\mbf A,\mbf b).$$ By $\ran(\mbf W_k)=\mcalk_k(\mbf A,\mbf b)$, we have \begin{align*} \mbf W_k^\top(\mbf b-\mbf A\mbf x_k^{\rm G})=\mbf W_k^\top(\mbf b-\mbf A\mbf W_k\mbf y_k^{\rm G})=\mbf W_k^\top(\mbf b-\mbf W_{k+1}\mbf H_{k+1,k}\mbf y_k^{\rm G})=\gamma_1\mbf e_1-\mbf H_k\mbf y_k^{\rm G}=\mbf 0.\end{align*} Greif and Varah \cite{greif2009itera} described a conjugate-gradient-type method (here we call it S$^2$CG) for nonsingular skew-symmetric linear systems.  S$^2$CG only computes the even iterates $\mbf x_{2j}^{\rm G}$ and returns $\mbf A^{-1}\mbf b$ in exact arithmetic. For singular skew-symmetric $\mbf A$, if $\mbf b\in\ran(\mbf A)$, then by Theorem \ref{oddeven}, $\ell$ is even, and thus S$^2$CG still applies; if $\mbf b\notin\ran(\mbf A)$, then $\ell$ is odd, and thus S$^2$CG does not apply. 

Now we consider Craig's method (CRAIG), which was originally derived as a form of CG applied to the normal equations of the second kind $$\mbf A\mbf A^\top\mbf y=\mbf b,\quad \mbf x=\mbf A^\top\mbf y.$$ Paige \cite{paige1974bidia} provided a description of CRAIG based on Golub--Kahan bidiagonalization (see section 2.2). The CRAIG points $\mbf x_j^{\rm CRAIG}$ satisfy $$ \mbf x_j^{\rm CRAIG}=\mbf V_j\mbf t_j= \mbf x_{j-1}^{\rm CRAIG}+\tau_j\mbf v_j,$$ where $\mbf t_j:=\bem\tau_1 & \tau_2 & \cdots & \tau_j\eem^\top$ satisfies $\mbf B_j\mbf t_j=\beta_1\mbf e_1$, and the components of $\mbf t_j$ can be computed recursively from $$\tau_1=\beta_1/\alpha_1,\quad \tau_j=-\beta_j\tau_{j-1}/\alpha_j\quad \text{for}\quad j\geq 2.$$ If $\mbf A$ is nonsingular, then CRAIG converges to the solution $\mbf A^{-1}\mbf b$. If $\mbf b\notin\ran(\mbf A)$, then CRAIG breaks down due to $\alpha_{k_0+1}=0$. Kammerer and Nashed \cite{kammerer1972conve} showed that in the singular consistent case, CG applied to the normal equations of the second kind converges to the solution $(\mbf A\mbf A^\top)^\dag\mbf b$. By using $ \mbf A^\top(\mbf A\mbf A^\top)^\dag\mbf b=\mbf A^\dag\mbf b$, we know that in the singular consistent case, CRAIG converges to the pseudoinverse solution of $\bf Ax=b$. 

For nonsingular skew-symmetric $\mbf A$, by using the zero structure in $\mbf H_k$ and taking $k=2j$, Greif et al. \cite[section 7]{greif2016numer} showed that $\mbf x_{2j}^{\rm G}=\mbf x_j^{\rm CRAIG}$. This also implies that $\mbf x_{2j}^{\rm G}\in\ran(\mbf A)$ since $\mbf x_j^{\rm CRAIG}=\mbf V_j\mbf t_j$ and $\ran(\mbf V_j)=\mcalk_j(\mbf A^\top\mbf A,\mbf A^\top\mbf b)=\mcalk_j(\mbf A^2,\mbf A\mbf b)$. Based on the same argument in \cite[section 7]{greif2016numer} and $\ran(\mbf A)=\nul(\mbf A)^\bot$, it is straightforward to prove the following.

\begin{proposition}\label{prop2}
Assume that $\mbf A$ is a singular skew-symmetric matrix, and that $\mbf b\in\ran(\mbf A)$. Let $\mbf x_j^{\rm G}$ and $\mbf x_j^{\rm CRAIG}$ be the $j$th iterates of {\rm S$^2$CG} and {\rm CRAIG} for $\mbf A\mbf x=\mbf b$, respectively. For each $j$ with $1\leq j\leq \ell/2$, we have $\mbf x_{2j}^{\rm G}=\mbf x_j^{\rm CRAIG}$. Moreover, {\rm S$^2$CG} returns the pseudoinverse solution $\mbf A^\dag\mbf b$. \end{proposition}

\subsection{S$^2$MR and LSQR}

Greif and Varah \cite{greif2009itera} described a minimum residual method (here we call it S$^2$MR) for a nonsingular skew-symmetric system. After $k$ iterations of Lanczos tridiagonalization, S$^2$MR seeks $$\mbf x_k^{\rm M}=\argmin_{\mbf x\in\mcalk_k(\mbf A,\mbf b)}\|\mbf b-\mbf A\mbf x\|_2.$$ By $\ran(\mbf W_k)=\mcalk_k(\mbf A,\mbf b)$, we have \begin{align*} \min_{\mbf x\in\mcalk_k(\mbf A,\mbf b)}\|\mbf b-\mbf A\mbf x\|_2 &=\min_{\mbf y\in\mbbr^k}\|\gamma_1\mbf w_1-\mbf A\mbf W_k\mbf y \|_2=\min_{\mbf y\in\mbbr^k}\|\gamma_1\mbf e_1-\mbf H_{k+1,k}\mbf y \|_2.\end{align*} The minimum residual solution is given by  $\mbf x_k^{\rm M}=\mbf W_k\mbf y_k^{\rm M}$ with $$\mbf y_k^{\rm M}=\argmin_{\mbf y\in\mbbr^k}\|\gamma_1\mbf e_1-\mbf H_{k+1,k}\mbf y\|_2.$$ Let the reduced QR factorization of $\mbf H_{k+1,k}$ be given by $$\mbf H_{k+1,k}=\mbf Q_{k+1,k}\mbf R_k.$$ We have $\mbf y_k^{\rm M}=\mbf R_k^{-1}\mbf Q_{k+1,k}^\top\gamma_1\mbf e_1$. If $\mbf z_k=\mbf Q_{k+1,k}^\top\gamma_1\mbf e_1$ and $\mbf P_k=\mbf W_k\mbf R_k^{-1}$, then $\mbf x_k^{\rm M}=\mbf P_k\mbf z_k$. The columns of $\mbf P_k$ and the entries of $\mbf z_k$ can be obtained by short recurrence relations; see \cite[section 4]{greif2009itera} for details.

Based on Golub--Kahan bidiagonalization, Paige and Saunders \cite{paige1982lsqr} proposed an iterative method named LSQR for solving linear systems and least-squares problems. LSQR is equivalent to CG applied to the normal equations $\mbf A^\top\mbf A\mbf x=\mbf A^\top\mbf b$. The LSQR points $\mbf x_j^{\rm LSQR}$ are characterized by \begin{align*}	\mbf x_j^{\rm LSQR}&=\argmin_{\mbf x\in\mcalk_j(\mbf A^\top\mbf A, \mbf A^\top\mbf b)}\|\mbf b-\mbf A\mbf x\|_2.\end{align*} We have $\mbf x_j^{\rm LSQR}=\mbf V_j\mbf y_j^{\rm LSQR}$ with $$\mbf y_j^{\rm LSQR}=\argmin_{\mbf y\in\mbbr^j}\|\beta_1\mbf e_1-\mbf B_{j+1,j}\mbf y\|_2.$$
For all types of linear systems (consistent or inconsistent, full column rank or rank-deficient), LSQR gives the pseudoinverse solution $\mbf A^\dag\mbf b$; see, e.g., \cite{fong2011lsmr}.

For nonsingular skew-symmetric $\mbf A$, by using the zero structure in $\mbf H_{k+1,k}$ and taking $k = 2j$, Greif et al. \cite[section 8]{greif2016numer} showed that $\mbf x_{2j}^{\rm M} =\mbf x_{2j+1}^{\rm M}=\mbf x_j^{\rm LSQR}$. This also implies that $\mbf x_{2j}^{\rm M}, \mbf x_{2j+1}^{\rm M}\in\ran(\mbf A)$ since $\mbf x_j^{\rm LSQR}=\mbf V_j\mbf y_j^{\rm LAQR}$ and $\ran(\mbf V_j)=\mcalk_j(\mbf A^\top\mbf A,\mbf A^\top\mbf b)=\mcalk_j(\mbf A^2,\mbf A\mbf b)$. Based on the same argument in \cite[section 8]{greif2016numer} and $\ran(\mbf A)=\nul(\mbf A)^\bot$, it is straightforward to prove the following.

\begin{proposition}\label{prop3}
Assume that $\mbf A$ is a singular skew-symmetric matrix. Let $\mbf x_j^{\rm M}$ and $\mbf x_j^{\rm LSQR}$ be the $j$th iterates of {\rm S$^2$MR} and {\rm LSQR} for $\mbf A\mbf x=\mbf b$, respectively. For each $j$ with $\mbf x_j^{\rm LSQR}\neq\mbf A^\dag\mbf b$, i.e., {\rm LSQR} does not converge at the $j$th iteration, we have $\mbf x_{2j}^{\rm M}=\mbf x_{2j+1}^{\rm M}=\mbf x_j^{\rm LSQR}$. 
Whether $\bf Ax=b$ is consistent or not, {\rm S$^2$MR} always returns the pseudoinverse solution $\mbf A^\dag\mbf b$. 
\end{proposition}
	
When applied to a singular symmetric linear system, MINRES gives the pseudoinverse solution if the system is consistent (see \cite[Theorem 2.4]{brown1997gmres} and \cite[Theorem 3.1]{choi2011minre}), and only gives a least-squares solution but not necessarily the pseudoinverse solution if the system is inconsistent (see \cite[Example 3.1 and Theorem 3.2]{choi2011minre}). To find the pseudoinverse solution to a singular inconsistent symmetric linear system, Choi et al. \cite {choi2011minre} proposed a MINRES-like method called MINRES-QLP. This method is also extended to solve complex symmetric, skew-symmetric, and skew-Hermitian systems \cite{choi2013minim}. Choi did not consider S$^2$MR. Unlike MINRES, our result shows that \rm S$^2$MR always returns the pseudoinverse solution $\mbf A^\dag\mbf b$, which was previously unnoticed.

\subsection{LSMR, LSLQ, and LNLQ}

LSMR \cite{fong2011lsmr}, LSLQ \cite{estrin2019lnlq}, LNLQ \cite{estrin2019lslq} are three other methods based on Golub--Kahan bidiagonalization. LSMR is equivalent to MINRES applied to the normal equations $\mbf A^\top\mbf A\mbf x=\mbf A^\top\mbf b$, and gives the pseudoinverse solution $\mbf A^\dag\mbf b$ of linear least-squares problems of any shape. LSLQ also solves linear least-squares problems of any shape, and is equivalent to SYMMLQ applied to the normal equations $\mbf A^\top\mbf A\mbf x=\mbf A^\top\mbf b$. LNLQ solves the linear least-norm problem $$\min\|\mbf x\|_2\quad  \text{ subject to }\quad \mbf A\mbf x=\mbf b,$$ and is equivalent to SYMMLQ applied to the normal equations of the second kind $$\mbf A\mbf A^\top\mbf y=\mbf b,\quad\mbf x=\mbf A^\top\mbf y.$$ All the three methods can be used to solve skew-symmetric linear systems. We summarize the convergence of {\rm S$^2$CG}, {\rm CRAIG}, {\rm S$^2$MR}, {\rm LSQR}, {\rm LSMR}, {\rm LSLQ}, and {\rm LNLQ} to the pseudoinverse solution $\mbf A^\dag\mbf b$ for all types of skew-symmetric linear systems in Table \ref{t1}. 
\begin{table}[htp]
\caption{Summary of the convergence of {\rm S$^2$CG},  {\rm CRAIG}, {\rm S$^2$MR}, {\rm LSQR}, {\rm LSMR}, {\rm LSLQ}, and {\rm LNLQ} to the pseudoinverse solution $\mbf A^\dag\mbf b$ for all types of skew-symmetric linear systems. {\rm Y} means the algorithm is convergent and {\rm N} means not.}  \label{t1}
\begin{center} 
\begin{tabular}{c|c|c|c|c|c|c|c} \toprule
Types  &  {\rm S$^2$CG} & {\rm CRAIG} & {\rm S$^2$MR} & {\rm LSQR} & {\rm LSMR} & {\rm LSLQ} & {\rm LNLQ} \\ 
\hline
singular consistent & Y & Y  &  Y  & Y  & Y & Y  & Y  \\ 
singular inconsistent & N & N  &  Y  & Y  & Y & Y  & N   \\
 nonsingular & Y & Y  &  Y  & Y  & Y & Y  & Y  \\ \bottomrule
\end{tabular}
\end{center}
\end{table}

For skew-symmetric linear systems, Greif et al. \cite[section 9]{greif2016numer} proposed a method based on Lanczos triangulation, and showed that it is numerically equivalent to LSMR by using the relation between Lanczos tridiagonalization and Golub--Kahan bidiagonalization for skew-symmetric matrices (see section 2.2). Based on the same argument in \cite[section 9]{greif2016numer}, we can propose two methods that are numerically equivalent to LSLQ and LNLQ, respectively. We omit the details, which will be given in Fan's master thesis \cite{fan2024itera}.

\section{Krylov subspace methods for shifted skew-symmetric linear systmes}

In this section, we mainly consider three Krylov subspaces methods (i.e., S$^3$CG, S$^3$MR, and S$^3$LQ) for shifted skew-symmetric linear systems. They are all based on Lanczos tridiagonalization, and correspond to CG, MINRES, and SYMMLQ for solving symmetric linear systems, respectively.

Next, we give the formulas established in Lanczos tridiagonalization for shifted skew-symmetric matrices. Assume that Lanczos tridiagonalization (Algorithm 1) applied to a skew-symmetric matrix $\mbf S$ and a nonzero vector $\mbf b$ stops in $\ell$ steps. We have $\gamma_k>0$ for $k=1,2\ldots,\ell$, $\gamma_{\ell+1}=0$, $\mbf w_{\ell+1}=\mbf 0$, and $$\mbf W_\ell^\top\mbf W_\ell=\mbf I_\ell,\qquad\mbf S\mbf W_\ell=\mbf W_\ell\mbf H_\ell,\qquad \mbf H_\ell=\mbf W_\ell^\top\mbf S\mbf W_\ell=\bem 0 & -\gamma_2 &&\\ \gamma_2 & 0 & \ddots &  \\ &\ddots &\ddots &-\gamma_\ell  \\ &&\gamma_\ell & 0 \eem.$$ 
 For the shifted skew-symmetric matrix $\mbf A=\alpha\mbf I+\mbf S$, we have $$\mbf A\mbf W_\ell=\alpha\mbf W_\ell+\mbf S\mbf W_\ell=\alpha\mbf W_\ell+\mbf W_\ell\mbf H_\ell=\mbf W_\ell\mbf T_\ell,\qquad \mbf T_\ell:=\alpha\mbf I_\ell+ \mbf H_\ell=\bem \alpha & -\gamma_2 &&\\ \gamma_2 & \alpha & \ddots &  \\ &\ddots &\ddots &-\gamma_\ell  \\ &&\gamma_\ell & \alpha \eem.$$ In the rest of this section, for each $k$ with $1\leq k\leq\ell-1$, let $\mbf T_k$, $\mbf T_{k+1,k}$, and $\mbf T_{k,k+1}$ denote the $k\times k$, $(k+1)\times k$,  and $k\times(k+1)$ upper-left submatrices of $\mbf T_\ell$, respectively. 
 
\subsection{S$^3$CG}

After $k$ iterations of Lanczos tridiagonalization, the Galerkin condition for solving the shifted skew-symmetric system $\bf Ax=b$ is \begin{equation}\label{galer}\mbf x_k^{\rm G}=\mbf W_k\mbf y_k^{\rm G},\qquad \mbf b-\mbf A\mbf x_k^{\rm G} \perp \mcalk_k(\mbf A,\mbf b).\end{equation} By $\ran(\mbf W_k)=\mcalk_k(\mbf A,\mbf b)$, we have \begin{align*} \mbf W_k^\top(\mbf b-\mbf A\mbf x_k^{\rm G})=\mbf W_k^\top(\mbf b-\mbf A\mbf W_k\mbf y_k^{\rm G})=\mbf W_k^\top(\mbf b-\mbf W_{k+1}\mbf T_{k+1,k}\mbf y_k^{\rm G})=\gamma_1\mbf e_1-\mbf T_k\mbf y_k^{\rm G}=\mbf 0.\end{align*} Since $\mbf T_k$ is nonsingular for each $k$, following the derivation of CG for symmetric positive definite systems (see, for example, \cite[Chapter 6]{demmel1997appli}), it is possible to devise a conjugate-gradient-type method for the shifted skew-symmetric system $\bf Ax=b$. We call the resulting method S$^3$CG and describe its implementation in Algorithm 4. 

\begin{center}
\begin{tabular*}{170mm}{l}
\toprule {\bf Algorithm 4}: S$^3$CG for $\bf Ax=b$ with shifted skew-symmetric ${\bf A}$ and nonzero $\mbf b$
\\ \hline \noalign{\smallskip}
\qquad Set $\mbf x_0^{\rm G}=\mbf 0$, $\mbf r_0^{\rm G}=\mbf b$ and $\mbf p_0^{\rm G}=\mbf r_0^{\rm G};$\\ \noalign{\smallskip}
\qquad {\bf for} $k=1,2,\ldots,$ {\bf do} until convergence:\\ \noalign{\smallskip}
\qquad\qquad $\alpha_k^{\rm G}=\dsp
\frac{(\mbf r_{k-1}^{\rm G})^\top\mbf r_{k-1}^{\rm G}}{(\mbf p_{k-1}^{\rm G})^\top\mbf A\mbf p_{k-1}^{\rm G}}$;\\ \noalign{\smallskip}
\qquad\qquad $\mbf x_k^{\rm G}=\mbf x_{k-1}^{\rm G}+\alpha_k^{\rm G}\mbf p_{k-1}^{\rm G}$; \\ \noalign{\smallskip}
\qquad\qquad $\mbf r_k^{\rm G}=\mbf r_{k-1}^{\rm G}-\alpha_k^{\rm G}\mbf A\mbf p_{k-1}^{\rm G}$;\\ \noalign{\smallskip}
\qquad\qquad $\beta_k^{\rm G}=-\dsp
\frac{(\mbf r_k^{\rm G})^\top\mbf r_k^{\rm G}}{(\mbf r_{k-1}^{\rm G})^\top\mbf r_{k-1}^{\rm G}}$;\\ \noalign{\smallskip}
\qquad\qquad $\mbf p_k^{\rm G}=\mbf r_k^{\rm G}+\beta_k^{\rm G}\mbf p_{k-1}^{\rm G}$;\\ \noalign{\smallskip}
\qquad {\bf end}\\
\bottomrule
\end{tabular*}
\end{center} 

The matrix $\mbf A=\alpha\mbf I+\mbf S$ can be regarded as a positive definite (see \cite{guducu2022nonhe}; if $\alpha<0$, use $-\mbf A$ instead) matrix. Applying Widlund's method \cite{widlund1978lancz} (also called the CGW method; see \cite[section 9.6]{saad2003itera}) to shifted skew-symmetric linear systems can yield a method that is mathematically equivalent to S$^3$CG. The approach here has the advantage that it emphasizes how S$^3$CG is related to CG for symmetric positive definite systems. The only difference is that the minus sign of $\beta_k^{\rm G}$ compared with CG. As in \cite[Chapter 6]{demmel1997appli} and \cite[Lecture 38]{trefethen1997numer}, we can prove the following results. 

\begin{proposition} Let {\rm S$^3$CG} be applied to a shifted skew-symmetric matrix problem $\bf Ax=b$. As long as the algorithm has not yet converged $($i.e., $\mbf r_{k-1}^{\rm G}\neq\mbf 0$$)$, it proceeds without divisions by zero, and we have the following identities of subspaces: \begin{align*}\mcalk_k({\bf A},{\bf b})&=\spa\{{\bf x}^{\rm G}_1,{\bf x}^{\rm G}_2,\cdots,{\bf x}^{\rm G}_k\}\\ &=\spa\{{\bf p}_0^{\rm G},{\bf p}_1^{\rm G},\cdots,{\bf p}_{k-1}^{\rm G}\}\\&=\spa\{{\bf r}_0^{\rm G},{\bf r}_1^{\rm G},\cdots,{\bf r}_{k-1}^{\rm G}\}.\end{align*} Moreover, the residuals are mutually orthogonal, $$(\mbf r_i^{\rm G})^\top\mbf r_k^{\rm G}=0\qquad (i\neq k),$$ and the search directions are ``semi-conjugate'', $$({\bf p}_i^{\rm G})^\top{\bf Ap}_k^{\rm G}=0\qquad (i<k).$$
\end{proposition}

\begin{remark} Note that the search directions $\mbf p_k$ in {\rm CG} for positive definite systems are conjugate, i.e, ${\bf p}_i^\top{\bf Ap}_k=0$ for all $i\neq k$. Here, we call the search directions $\mbf p_k^{\rm G}$ in {\rm S$^3$CG} for shifted skew-symmetric systems  ``semi-conjugate'' because $({\bf p}_i^{\rm G})^\top{\bf Ap}_k^{\rm G}=0$ holds only for $i<k$.
\end{remark}

The optimality properties and error estimates for Widlund's method have been given in \cite{eisenstat1983note,hageman1980equiv,szyld1993varia}. Applying these results to shifted skew-symmetric linear systems, we can obtain the optimality properties and error estimates for S$^3$CG, which will be given below. The even and odd subsequences $\{\mbf x_{2k}^{\rm G}\}$ and $\{\mbf x_{2k+1}^{\rm G}\}$ of {\rm S$^3$CG} satisfy the optimality properties \beq\label{x2kG}\|\mbf x_{2k}^{\rm G}-\mbf A^{-1}\mbf b\|_2=\min_{\mbf x\in\mbf A^\top\mcalk_{2k}(\mbf A,\mbf b)}\|\mbf x-\mbf A^{-1}\mbf b\|_2,\eeq and $$\|\mbf x_{2k+1}^{\rm G}-\mbf A^{-1}\mbf b\|_2=\min_{\mbf x\in\mbf b/\alpha+\mbf A^\top\mcalk_{2k+1}(\mbf A,\mbf b)}\|\mbf x-\mbf A^{-1}\mbf b\|_2.$$	The eigenvalues of $\mbf S$ are purely imaginary. Let $\rmi[-\beta, \beta]$ for some $\beta>0$ be the smallest interval that contains these eigenvalues. The optimality property of S$^3$CG leads to an error bound of the form $$\frac{\|\mbf x_{2k}^{\rm G}-\mbf A^{-1}\mbf b\|_2}{\|\mbf A^{-1}\mbf b\|_2}\leq 2\l(\frac{\sqrt{1+(\beta/\alpha)^2}-1}{\sqrt{1+(\beta/\alpha)^2}+1}\r)^k.$$ 
The same bound holds for the sequence $\|\mbf x_{2k+1}^{\rm G}-\mbf A^{-1}\mbf b\|_2/\|\mbf x_1^{\rm G}-\mbf A^{-1}\mbf b\|_2$. The bound indicates that a ``fast" convergence of S$^3$CG can be expected when $\beta/|\alpha|> 0$ is ``small".

\subsubsection{Relation to CRAIG}

In this subsection we compare S$^3$CG with CRAIG for shifted skew-symmetric linear systems. We will show that even iterations of S$^3$CG are CRAIG iterations. The following lemma is important for our discussion, which is an extension of Lemma 1 of \cite {eisenstat2015equiv} to shifted skew-symmetric matrices.

\begin{lemma}\label{orth} Let $\mbf A=\alpha\mbf I+\mbf S$ be a shifted skew-symmetric matrix. The subspaces $\mbf A^\top\mcalk_k(\mbf S^2,\mbf b)$ and $\mbf A^\top\mcalk_k(\mbf S^2,\mbf S\mbf b)$ are orthogonal and the solution $\mbf A^{-1}\mbf b$ is orthogonal to $\mbf A^\top\mcalk_k(\mbf S^2,\mbf S\mbf b)$. 	 
\end{lemma}
\proof It is sufficient to show that both $\mbf A^\top(\mbf S^2)^i\mbf b$ and $\mbf A^{-1}\mbf b$ are orthogonal to $\mbf A^\top(\mbf S^2)^j\mbf S\mbf b$ for any $0\leq i,j\leq k-1$. Note that  $$(\mbf A^{-1}\mbf b)^\top\mbf A^\top(\mbf S^2)^j\mbf S\mbf b=\mbf b^\top\mbf S^{2j+1}\mbf b=(-1)^j(\mbf S^j\mbf b)^\top\mbf S(\mbf S^j\mbf b)=0$$ and \begin{align*}(\mbf A^\top(\mbf S^2)^i\mbf b)^\top\mbf A^\top(\mbf S^2)^j\mbf S\mbf b & =\mbf b^\top\mbf S^{2i}\mbf A\mbf A^\top\mbf S^{2j+1}\mbf b=\mbf b^\top\mbf S^{2i}(\alpha^2\mbf I-\mbf S^2)\mbf S^{2j+1}\mbf b\\ &= \alpha^2\mbf b^\top\mbf S^{2i+2j+1}\mbf b-\mbf b^\top\mbf S^{2i+2j+3}\mbf b\\ &= \alpha^2(-1)^{i+j}(\mbf S^{i+j}\mbf b)^\top\mbf S(\mbf S^{i+j}\mbf b)-(-1)^{i+j+1}(\mbf S^{i+j+1}\mbf b)^\top\mbf S(\mbf S^{i+j+1}\mbf b)\\ &=0\end{align*} for any $0\leq i, j\leq k-1$. This completes the proof.
\endproof

Now we give the main theoretical result of this subsection. 

\begin{theorem}\label{theo7} Let $\mbf A=\alpha\mbf I+\mbf S$ be a shifted skew-symmetric matrix. Let $\mbf x_k^{\rm G}$ and $\mbf x_k^{\rm CRAIG}$ be the $k$th iterates of {\rm S$^3$CG} and {\rm CRAIG} for $\mbf A\mbf x=\mbf b$, respectively.  Then we have
$$\mbf x_{2k}^{\rm G}=\mbf x_k^{\rm CRAIG}.$$
\end{theorem}
\proof The CRAIG points $\mbf x_k^{\rm CRAIG}$ are characterized by (see, for example, \cite[Chapter 7]{greenbaum1997itera}) \begin{align*}\mbf x_k^{\rm CRAIG}=\argmin_{\mbf x\in\mbf A^\top\mcalk_k(\mbf A\mbf A^\top, \mbf b)}\|\mbf x-\mbf A^{-1}\mbf b\|_2^2=\argmin_{\mbf x\in\mbf A^\top\mcalk_k(\mbf S^2, \mbf b)}\|\mbf x-\mbf A^{-1}\mbf b\|_2^2.\end{align*} The last equality is a direct result of $\mbf A\mbf A^\top=\alpha^2\mbf I-\mbf S^2$ and $\mcalk_k(\alpha^2\mbf I-\mbf S^2,\mbf b)=\mcalk_k(\mbf S^2,\mbf b)$. The S$^3$CG points $\mbf x_{2k}^{\rm G}$ are characterized by (see (\ref{x2kG})) $$\mbf x_{2k}^{\rm G}=\argmin_{\mbf x\in\mbf A^\top\mcalk_{2k}(\mbf A,\mbf b)}\|\mbf x-\mbf A^{-1}\mbf b\|_2^2=\argmin_{\mbf x\in\mbf A^\top\mcalk_{2k}(\mbf S,\mbf b)}\|\mbf x-\mbf A^{-1}\mbf b\|_2^2.$$ The last equality is a direct result of $\mcalk_{2k}(\alpha\mbf I+\mbf S,\mbf b)=\mcalk_{2k}(\mbf S,\mbf b)$. For any $\mbf x\in\mbf A^\top\mcalk_{2k}(\mbf S,\mbf b)$, by $$\mbf A^\top\mcalk_{2k}(\mbf S,\mbf b)=\mbf A^\top\mcalk_k(\mbf S^2,\mbf b) + \mbf A^\top\mcalk_k(\mbf S^2,\mbf S\mbf b),$$ we can write $\mbf x=\mbf y+\mbf z$, where $\mbf y\in\mbf A^\top\mcalk_k(\mbf S^2,\mbf b)$ and $\mbf z\in\mbf A^\top\mcalk_k(\mbf S^2,\mbf S\mbf b)$. By Lemma \ref{orth} and the Pythagorean Theorem, we have \begin{align*}\min_{\mbf x\in\mbf A^\top\mcalk_{2k}(\mbf S,\mbf b)}\|\mbf x-\mbf A^{-1}\mbf b\|_2^2 &=\min_{\mbf y\in\mbf A^\top\mcalk_k(\mbf S^2,\mbf b)}\|\mbf y-\mbf A^{-1}\mbf b\|_2^2+\min_{\mbf z\in\mbf A^\top\mcalk_k(\mbf S^2,\mbf S\mbf b)}\|\mbf z\|_2^2\\ &=\min_{\mbf y\in\mbf A^\top\mcalk_k(\mbf S^2,\mbf b)}\|\mbf y-\mbf A^{-1}\mbf b\|_2^2 ,\end{align*} and hence $\mbf x_{2k}^{\rm G}=\mbf x_k^{\rm CRAIG}$.  \endproof
 
\subsection{S$^3$MR}

After $k$ iterations of Lanczos tridiagonalization, the minimum residual condition is $$\mbf x_k^{\rm M}=\argmin_{\mbf x\in\mcalk_k(\mbf A,\mbf b)}\|\mbf b-\mbf A\mbf x\|_2.$$ We have $\mbf x_k^{\rm M}=\mbf W_k\mbf y_k^{\rm M}$, where $\mbf y_k^{\rm M}\in\mbbr^k$ solves $$\min_{\mbf y\in\mbbr^k}\|\gamma_1\mbf e_1-\mbf T_{k+1,k}\mbf y\|_2.$$ Similar to MINRES for symmetric systems, an implementation using short recurrences has been given in \cite{idema2007minim,jiang2007algor,idema2023compa}. For completeness, we also provide the implementation of S$^3$MR in Algorithm 5. We note that applying Rapoport's method \cite{rapoport1978nonli} to shifted skew-symmetric linear systems can yield a method that is mathematically equivalent to S$^3$MR.

\begin{center}
\begin{tabular*}{170mm}{l}
\toprule {\bf Algorithm 5}: S$^3$MR for $\bf Ax=b$ with shifted skew-symmetric ${\bf A}=\alpha\mbf I+\mbf S$ and nonzero $\mbf b$ 
\\ \hline \noalign{\smallskip}
\qquad Set $\mbf x_0^{\rm M}=\mbf 0$, $\wt\delta_1=\alpha$, $c_0=1$,  $\gamma_1=\|\mbf b\|_2$, $\wt\psi_1=\gamma_1$, $\mbf w_0=\mbf 0$, and $\mbf w_1=\mbf b/\gamma_1$;\\ \noalign{\smallskip}
\qquad {\bf for} $k=1,2,\ldots,$ {\bf do} until convergence:\\ \noalign{\smallskip}
\qquad\qquad $\gamma_{k+1}\mbf w_{k+1}:=\mbf S\mbf w_k+\gamma_k\mbf w_{k-1}$;\\ \noalign{\smallskip}
\qquad\qquad $\delta_k=\sqrt{\wt\delta_k^2+\gamma_{k+1}^2}$, $c_k=\wt\delta_k/\delta_k$, $s_k=\gamma_{k+1}/\delta_k$;\\ \noalign{\smallskip}
\qquad\qquad  $\wt\delta_{k+1}=\alpha c_k+\gamma_{k+1}c_{k-1}s_k$; \\ \noalign{\smallskip}
\qquad\qquad  $\psi_k=c_k\wt\psi_k$, $ \wt\psi_{k+1}=-s_k\wt\psi_k$; \\ \noalign{\smallskip}
\qquad\qquad {\bf if} $k\leq 2$ {\bf then} \\ \noalign{\smallskip}
\qquad\qquad\qquad $\mbf p_k=\mbf w_k/\delta_k$;  \\ \noalign{\smallskip}
\qquad\qquad {\bf else} \\ \noalign{\smallskip}
\qquad\qquad \qquad $\mbf p_k=(\mbf w_k+\gamma_ks_{k-2}\mbf p_{k-2})/\delta_k$;\\ \noalign{\smallskip}
\qquad\qquad {\bf end}\\ \noalign{\smallskip}
\qquad\qquad $\mbf x_k^{\rm M} =  \mbf x_{k-1}^{\rm M}+ \psi_k\mbf p_k$;\\ \noalign{\smallskip}
\qquad {\bf end}\\
\bottomrule
\end{tabular*}
\end{center}

The residual $\mbf r_k^{\rm M}=:\mbf b-\mbf A\mbf x_k^{\rm M}$ of S$^3$MR satisfies (see \cite{jiang2007algor})  
 $$\|\mbf r_k^{\rm M}\|_2=|\wt\psi_{k+1}|=\cdots=|s_ks_{k-1}\cdots s_1|\gamma_1.$$ Since $\mbf T_k$ is nonsingular for each $k$, S$^3$MR does not stagnate (see, for example, \cite{meurant2014neces} and \cite[section 2]{du2017any}), i.e., $\|\mbf r_k^{\rm M}\|_2$ is strictly decreasing.
Let $\rmi[-\beta, \beta]$ for some $\beta>0$ be the smallest interval that contains the eigenvalues of the skew-symmetric matrix $\mbf S$.  The optimality property of S$^3$MR leads to a convergence bound (see \cite[Theorem 4.3]{szyld1993varia} and \cite[section 4]{jiang2007algor}) $$\frac{\|\mbf r_k^{\rm M}\|_2}{\|\mbf b\|_2}\leq 2\l(\frac{\beta/|\alpha|}{\sqrt{1+(\beta/\alpha)^2}+1}\r)^k.$$  
As for S$^3$CG, this bound for S$^3$MR indicates that the convergence is ``fast'' when $\beta/|\alpha|>0$ is ``small''.

\subsubsection{Relation to LSQR}

In this subsection we compare S$^3$MR with LSQR for shifted skew-symmetric linear systems. We note that Golub--Kahan bidiagonalization (Algorithm 2) for shifted skew-symmetric $\mbf A=\alpha\mbf I+\mbf S$ ($\alpha\neq0$) and nonzero $\mbf b$ must stop in $\ell_0=\lceil\ell/2\rceil$ steps with $\beta_{\ell_0+1}=0$ because $\mbf b\in\ran(\mbf A)$ always holds. An interesting theoretical result on Lanczos tridiagonalization and Golub--Kahan bidiagonalization is presented in the appendix.

\begin{proposition}\label{leq} Let $\mbf A=\alpha\mbf I+\mbf S$ and $\alpha\neq 0$. Let $\mbf x_k^{\rm M}$ and $\mbf x_k^{\rm LSQR}$ be the $k$th iterates of {\rm S$^3$MR} and {\rm LSQR} for $\mbf A\mbf x=\mbf b$, respectively. For each $k$ with $1\leq k\leq \ell_0-1$, it holds that $$\|\mbf b-\mbf A\mbf x_{2k}^{\rm M}\|_2\leq\|\mbf b-\mbf A\mbf x_k^{\rm LSQR}\|_2.$$ Moreover, we have $\mbf x_\ell^{\rm M}=\mbf x_{\ell_0}^{\rm LSQR}=\mbf A^{-1}\mbf b$.
\end{proposition}

\proof
For LSQR, we have $$\mbf x_k^{\rm LSQR}=\argmin_{\mbf x\in\mcalk_k(\mbf A^\top\mbf A,\mbf A^\top\mbf b)}\|\mbf b-\mbf A\mbf x\|_2=\argmin_{\mbf x\in\mcalk_k(\mbf S^2,\alpha\mbf b-\mbf S\mbf b)}\|\mbf b-\mbf A\mbf x\|_2.$$ For S$^3$MR, we have $$\mbf x_{2k}^{\rm M}=\argmin_{\mbf x\in\mcalk_{2k}(\mbf A,\mbf b)}\|\mbf b-\mbf A\mbf x\|_2=\argmin_{\mbf x\in\mcalk_{2k}(\mbf S,\mbf b)}\|\mbf b-\mbf A\mbf x\|_2.$$ It follows from $\mcalk_k(\mbf S^2,\alpha\mbf b-\mbf S\mbf b)\subset\mcalk_{2k}(\mbf S,\mbf b)$ that $\|\mbf b-\mbf A\mbf x_{2k}^{\rm M}\|_2\leq\|\mbf b-\mbf A\mbf x_k^{\rm LSQR}\|_2$. By $\gamma_{\ell+1}=0$ and $\beta_{\ell_0+1}=0$, it is straightforward to show that $\mbf x_\ell^{\rm M}=\mbf x_{\ell_0}^{\rm LSQR}=\mbf A^{-1}\mbf b$. 
\endproof

Proposition \ref{leq} shows that even residuals of S$^3$MR are not greater than residuals of LSQR. Our numerical experiments (see section 5) show that the strict inequality holds true, namely, that $\|\mbf b-\mbf A\mbf x_{2k}^{\rm M}\|_2<\|\mbf b-\mbf A\mbf x_k^{\rm LSQR}\|_2$.

\subsection{S$^3$LQ}
S$^3$LQ fills a gap in the family of iterative methods for shifted skew-symmetric linear systems based on Lanczos tridiagonalization. Similar to SYMMLQ for symmetric linear systems, we can defined the S$^3$LQ points to be 
\begin{align*}
 	\mbf x_k^{\rm L}:=\argmin_{\mbf x\in\mcalk_k(\mbf A,\mbf b)}\|\mbf x\|_2\quad \mbox{subject to}\quad \mbf b-\mbf A\mbf x \perp \mcalk_{k-1}(\mbf A,\mbf b).
 \end{align*} By Lanczos tridiagonalization, we have $\mbf x_k^{\rm L}=\mbf W_k\mbf y_k^{\rm L}$, where $\mbf y_k^{\rm L}\in\mbbr^k$ solves \begin{equation}\label{ycons}\min \|\mbf y\|_2\quad {\rm subject\ to}\quad \mbf T_{k-1,k}\mbf y=\gamma_1\mbf e_1.\end{equation} Similar to SYMMLQ, the error $\|\mbf x_k^{\rm L}-\mbf A^{-1}\mbf b\|_2$ is monotonically decreasing. In fact, we have the following optimality result.
 
 \begin{theorem}\label{theo9} For $k>1$, we have $$\mbf x_k^{\rm L}=\argmin_{\mbf x\in\mbf A^\top\mcalk_{k-1}(\mbf A,\mbf b)}\|\mbf x-\mbf A^{-1}\mbf b\|_2.$$ 
 \end{theorem}
 \proof From the theory of least squares problems we know that the optimal solution of $$\min_{\mbf x\in\mbf A^\top\mcalk_{k-1}(\mbf A,\mbf b)}\|\mbf x-\mbf A^{-1}\mbf b\|_2$$ is $$\mbf A^\top\mbf W_{k-1}(\mbf A^\top\mbf W_{k-1})^\dag\mbf A^{-1}\mbf b.$$ Similarly, the solution of (\ref{ycons}) is $\mbf y_k^{\rm L}=\mbf T_{k-1,k}^\dag\gamma_1\mbf e_1$.  By $\mbf A^\top\mbf W_\ell=\mbf W_\ell\mbf T_\ell^\top$, we have \begin{align*}\mbf A^\top\mbf W_{k-1}(\mbf A^\top\mbf W_{k-1})^\dag\mbf A^{-1}\mbf b&=\mbf W_k\mbf T_{k-1,k}^\top(\mbf T_{k-1,k}\mbf W_k^\top\mbf W_k\mbf T_{k-1,k}^\top)^{-1}\mbf W_{k-1}^\top\mbf A\mbf A^{-1}\mbf b\\ & = \mbf W_k\mbf T_{k-1,k}^\top(\mbf T_{k-1,k}\mbf T_{k-1,k}^\top)^{-1}\mbf W_{k-1}^\top\mbf b\\ &=\mbf W_k\mbf T_{k-1,k}^\dag\gamma_1\mbf e_1=\mbf W_k\mbf y_k^{\rm L}=\mbf x_k^{\rm L}.\end{align*} This completes the proof.
 \endproof

We now discuss the implementation of S$^3$LQ. The vector $\mbf y_k^{\rm L}$ can be obtained via the LQ factorization \begin{equation}\label{lq}\mbf T_{k-1,k}\mbf Q^\top_k=\bem\mbf L_{k-1} &\mbf 0\eem, \qquad \mbf L_{k-1}:=\bem \delta_1 &&&&\\ \lambda_2 &\delta_2 &&&\\ \eta_3 & \lambda_3 &\delta_3 &&\\ &\ddots &\ddots &\ddots &\\ && \eta_{k-1} & \lambda_{k-1}&\delta_{k-1}\eem,\end{equation} where $$\mbf Q_k^\top=\mbf G_{1,2}\mbf G_{2,3}\cdots\mbf G_{k-1,k}$$ is a product of Givens rotations, and $\mbf L_{k-1}$ is lower triangular. The $i$th Givens rotation $\mbf G_{i,i+1}$ can be represented as $$\mbf G_{i,i+1}=\bem \mbf I_{i-1} & & &\\ & c_i & -s_i &\\ & s_i & c_i & \\ &&& \mbf I_{k-1-i}\eem\in\mbbr^{k\times k},$$ which zeros out $-\gamma_{i+1}$. For the purpose of establishing recursion formulae, we define $\wt\delta_1:=\alpha$ and $c_0:=1$. We have $$\bem \wt\delta_i & -\gamma_{i+1}\\ \gamma_{i+1}c_{i-1} & \alpha \\ 0 & \gamma_{i+2} \eem\bem c_i & -s_i\\ s_i & c_i \eem = \bem \delta_i & 0\\ \lambda_{i+1} & \wt\delta_{i+1} \\ \eta_{i+2}& \gamma_{i+2}c_i\eem.$$ The $i$th Givens rotation zeros out $-\gamma_{i+1}$, and hence we have $$\delta_i:=\sqrt{\wt\delta_i^2+\gamma_{i+1}^2}, \qquad c_i:=\wt\delta_i/\delta_i,\qquad s_i=-\gamma_{i+1}/\delta_i,$$ and $$ \lambda_{i+1}= \gamma_{i+1}c_{i-1}c_i+\alpha s_i, \qquad \wt\delta_{i+1}=\alpha c_i-\gamma_{i+1}c_{i-1}s_i,\qquad \eta_{i+2}=\gamma_{i+2}s_i.$$ The following result shows that $\mbf L_{k-1}$ has only two nonzero  diagonals.

\begin{proposition}
	We have $\lambda_2=\lambda_3=\cdots=\lambda_{k-1}=0.$
\end{proposition}
\proof Direct calculations yield $$\lambda_2=\gamma_2c_0c_1+\alpha s_1=\gamma_2\alpha/\delta_1+\alpha(-\gamma_2/\delta_1)=0,$$ and
\begin{align*}\lambda_{i+1} &= \gamma_{i+1}c_{i-1}c_i+\alpha s_i  = \frac{\gamma_{i+1}}{\delta_i}(c_{i-1}\wt\delta_i-\alpha)\\ &= \frac{\gamma_{i+1}}{\delta_i}(c_{i-1}^2\alpha-\gamma_i c_{i-2}c_{i-1}s_{i-1} -\alpha)\\ &= \frac{\gamma_{i+1}}{\delta_i}(-s_{i-1}^2\alpha-\gamma_ic_{i-2}c_{i-1}s_{i-1} )\\ &= \frac{\gamma_{i+1}}{\delta_i}(-s_{i-1})(\gamma_ic_{i-2}c_{i-1}+\alpha s_{i-1})\\ &= \frac{\gamma_{i+1}}{\delta_i}(-s_{i-1})\lambda_i.\end{align*} This completes the proof. 
\endproof

Now we have
$$\mbf L_{k-1}=\bem \delta_1 &&&&\\ 0 &\delta_2 &&&\\ \eta_3 & 0 &\delta_3 &&\\ &\ddots &\ddots &\ddots &\\ && \eta_{k-1} & 0 &\delta_{k-1} \eem.$$ The factorization (\ref{lq}) allows us to rewrite the constraints of (\ref{ycons}) as $$\mbf L_{k-1}\mbf z_{k-1}=\gamma_1\mbf e_1$$ and hence the solution of (\ref{ycons}) is $$\mbf y_k^{\rm L}=\mbf Q_k^\top\bem\mbf z_{k-1}\\ 0\eem.$$ If $\mbf z_{k-1}:=\bem\xi_1 & \xi_2&\cdots & \xi_{k-1}\eem^\top$, then  \begin{equation}\label{xieven}\xi_1=\gamma_1/\delta_1, \qquad \xi_2=0,\qquad \xi_{k-1}=-\eta_{k-1}\xi_{k-3}/\delta_{k-1}\quad (k\geq 4).\end{equation} Obviously, it holds that $\mbf z_{k-1}=\bem\mbf z_{k-2}^\top & \xi_{k-1}\eem^\top$. Next we define \begin{align*}\mbf P_k:=\mbf W_k\mbf Q_k^\top =\bem \mbf p_1 & \mbf p_2 &\cdots & \mbf p_{k-1} & \wt{\mbf p}_k\eem.\end{align*} 
By using \begin{align*}
\mbf W_k=\bem\mbf W_{k-1} & \mbf w_k\eem,\quad  \mbf Q_k^\top=\bem\mbf Q_{k-1}^\top &\\& 1\eem\mbf G_{k-1,k},\end{align*} we have
\begin{align*}
\mbf P_k&=\bem\mbf W_{k-1}\mbf Q_{k-1}^\top & \mbf w_k\eem\mbf G_{k-1,k}=\bem\mbf P_{k-1}& \mbf w_k\eem\mbf G_{k-1,k}\\ &=\bem \mbf p_1 & \mbf p_2 &\cdots & \mbf p_{k-2} & \wt{\mbf p}_{k-1} & \mbf w_k \eem\bem\mbf I_{k-2}&&\\&c_{k-1} & -s_{k-1}\\& s_{k-1}& c_{k-1}\eem,\end{align*} which gives $$\mbf p_{k-1} = c_{k-1} \wt{\mbf p}_{k-1}+s_{k-1}\mbf w_k,\qquad \wt{\mbf p}_k = c_{k-1} \mbf w_k-s_{k-1}\wt{\mbf p}_{k-1}.$$ The approximate solution $\mbf x_k^{\rm L}$ may be updated efficiently by using the recursion for $k\geq 2$, $$\mbf x_k^{\rm L}=\mbf W_k\mbf y_k^{\rm L}=\mbf W_k\mbf Q_k^\top\bem\mbf z_{k-1}\\ 0\eem=\mbf P_k\bem\mbf z_{k-1}\\ 0\eem= \mbf P_{k-1}\bem\mbf z_{k-2}\\ 0\eem+ \xi_{k-1}\mbf p_{k-1}= \mbf x_{k-1}^{\rm L}+ \xi_{k-1}\mbf p_{k-1},$$  with $\mbf x_1^{\rm L}=\mbf 0$ and $\wt{\mbf p}_1=\mbf w_1$. Since $\mbf x_k^{\rm L}$ is updated along orthogonal directions, $\|\mbf x_k^{\rm L}\|_2$ is monotonically increasing. Algorithm 6 summarizes S$^3$LQ.

\begin{center}
\begin{tabular*}{170mm}{l}
\toprule {\bf Algorithm 6}: S$^3$LQ for $\bf Ax=b$ with shifted skew-symmetric ${\bf A}=\alpha\mbf I+\mbf S$ and nonzero $\mbf b$ 
\\ \hline \noalign{\smallskip}
\qquad Set $\mbf x_1^{\rm L}=\mbf 0$, $\wt\delta_1=\alpha$, $c_0=1$,  $\gamma_1=\|\mbf b\|_2$, $\mbf w_0=\mbf 0$, $\mbf w_1=\mbf b/\gamma_1$, and $\wt{\mbf p}_1=\mbf w_1$;\\ \noalign{\smallskip}
\qquad {\bf for} $k=1,2,\ldots,$ {\bf do} until convergence:\\ \noalign{\smallskip}
\qquad\qquad $\gamma_{k+1}\mbf w_{k+1}:=\mbf S\mbf w_k+\gamma_k\mbf w_{k-1}$;\\ \noalign{\smallskip}
\qquad\qquad $\delta_k=\sqrt{\wt\delta_k^2+\gamma_{k+1}^2}$, $c_k=\wt\delta_k/\delta_k$, $s_k=-\gamma_{k+1}/\delta_k$;\\ \noalign{\smallskip}
\qquad\qquad  $\wt\delta_{k+1}=\alpha c_k-\gamma_{k+1}c_{k-1}s_k$; \\ \noalign{\smallskip}
\qquad\qquad {\bf if} $k=1$ {\bf then} \\ \noalign{\smallskip}
\qquad\qquad\qquad $\xi_1=\gamma_1/\delta_1$;  \\ \noalign{\smallskip}
\qquad\qquad {\bf else if} $k=2$ {\bf then} \\ 
\qquad\qquad \qquad $\xi_2=0$;  \\ \noalign{\smallskip} 
\qquad\qquad {\bf else} \\ \noalign{\smallskip}
\qquad\qquad \qquad $\xi_k=-\gamma_ks_{k-2}\xi_{k-2}/\delta_k$;\\ \noalign{\smallskip}
\qquad\qquad {\bf end}\\ \noalign{\smallskip}
\qquad\qquad $\mbf p_k = c_k \wt{\mbf p}_k+s_k\mbf w_{k+1}$; \\ \noalign{\smallskip}
\qquad\qquad $\mbf x_{k+1}^{\rm L} =  \mbf x_k^{\rm L}+ \xi_k\mbf p_k$;\\ \noalign{\smallskip}
\qquad\qquad $\wt{\mbf p}_{k+1} = c_k \mbf w_{k+1}-s_k\wt{\mbf p}_k$\\ \noalign{\smallskip}
\qquad {\bf end}\\
\bottomrule
\end{tabular*}
\end{center}

The residual \begin{align*}\mbf r_k^{\rm L}&:=\mbf b-\mbf A\mbf x_k^{\rm L}=\mbf b-\mbf A\mbf W_k\mbf y_k^{\rm L} =\mbf W_k\gamma_1\mbf e_1-(\mbf W_k\mbf T_k+\gamma_{k+1}\mbf w_{k+1}\mbf e_k^\top)\mbf y_k^{\rm L}\\  &=\mbf W_k(\gamma_1\mbf e_1-\mbf T_k\mbf y_k^{\rm L})-\gamma_{k+1}\mbf w_{k+1}\mbf e_k^\top\mbf y_k^{\rm L}\\ &=-(\gamma_k\mbf e_{k-1}^\top\mbf y_k^{\rm L}+\alpha\mbf e_k^\top\mbf y_k^{\rm L})\mbf w_k-\gamma_{k+1}\mbf e_k^\top\mbf y_k^{\rm L}\mbf w_{k+1}.\end{align*} It follows from \begin{align*}\mbf e_{k-1}^\top\mbf y_k^{\rm L} &=\mbf e_{k-1}^\top\mbf G_{k-2,k-1}\mbf G_{k-1,k}\bem\mbf z_{k-1}\\  0\eem=(s_{k-2}\mbf e_{k-2}^\top+c_{k-2}c_{k-1}\mbf e_{k-1}^\top-c_{k-2}s_{k-1}\mbf e_k^\top)\bem\mbf z_{k-1}\\  0\eem\\ &=s_{k-2}\xi_{k-2}+c_{k-2}c_{k-1}\xi_{k-1}\end{align*} and \begin{align*} \mbf e_k^\top\mbf y_k^{\rm L}=\mbf e_k^\top\mbf G_{k-1,k}\bem\mbf z_{k-1}\\ 0\eem=(s_{k-1}\mbf e_{k-1}^\top+c_{k-1}\mbf e_k^\top)\bem\mbf z_{k-1}\\ 0\eem=s_{k-1}\xi_{k-1} \end{align*} that \begin{align*}\mbf r_k^{\rm L}&=-(\gamma_k s_{k-2}\xi_{k-2}+\gamma_kc_{k-2}c_{k-1}\xi_{k-1}+\alpha s_{k-1}\xi_{k-1})\mbf w_k-\gamma_{k+1}s_{k-1}\xi_{k-1}\mbf w_{k+1}\\ &=-(\gamma_ks_{k-2}\xi_{k-2}+\lambda_k\xi_{k-1})\mbf w_k-\gamma_{k+1}s_{k-1}\xi_{k-1}\mbf w_{k+1}\\ & = -\eta_k\xi_{k-2}\mbf w_k-\eta_{k+1}\xi_{k-1}\mbf w_{k+1}.\end{align*} By $\xi_{2j}=0$ (see (\ref{xieven})) we have for $j\geq 1$ $$\mbf r_{2j}^{\rm L}=\mbf r_{2j+1}^{\rm L}=-\eta_{2j+1}\xi_{2j-1}\mbf w_{2j+1}.$$ It follows that $$\|\mbf r_{2j}^{\rm L}\|_2=\|\mbf r_{2j+1}^{\rm L}\|_2=|\eta_{2j+1}\xi_{2j-1}|.$$

Now we consider the Galerkin approximation solution $\mbf x_k^{\rm G}$; see (\ref{galer}). Direct calculations yields the LQ factorization $$\mbf T_k\mbf Q_k^\top=\wt{\mbf L}_k=\bem \delta_1 &&&&\\ 0 &\delta_2 &&&\\ \eta_3 & 0 &\delta_3 &&\\ &\ddots &\ddots &\ddots &\\ && \eta_k & 0 &\wt\delta_k \eem,$$ where $\wt{\mbf L}_k$ differs from $\mbf L_k$ only in the $(k,k)$th element. This factorization allows us to rewrite $\mbf T_k\mbf y_k^{\rm G}=\gamma_1\mbf e_1$ as $$\wt{\mbf L}_k\wt{\mbf z}_k=\gamma_1\mbf e_1,\qquad \mbf y_k^{\rm G}=\mbf Q_k^\top\wt{\mbf z}_k.$$ By the structure of $\wt{\mbf L}_k$, i.e., $$\wt{\mbf L}_k=\bem \mbf L_{k-1} & \mbf 0\\ \eta_k\mbf e_{k-2}^\top &\wt\delta_k\eem,$$ and $\mbf L_{k-1}\mbf z_{k-1}=\gamma_1\mbf e_1$, $\mbf z_{k-1}=\bem\xi_1 & \xi_2 & \cdots &\xi_{k-1}\eem^\top$, we have \beq\label{xi0}\wt{\mbf z}_k=\bem\mbf z_{k-1}\\ \wt\xi_k\eem,\qquad \wt\xi_k=-\eta_k\xi_{k-2}/\wt\delta_k.\eeq Then, by using $\mbf W_k\mbf Q_k^\top =\bem \mbf p_1 & \mbf p_2 &\cdots & \mbf p_{k-1} & \wt{\mbf p}_k\eem$  we have $$\mbf x_k^{\rm G}=\mbf W_k\mbf y_k^{\rm G}=\mbf W_k\mbf Q_k^\top \bem\mbf z_{k-1}\\ \wt\xi_k\eem=\mbf W_k\mbf Q_k^\top\bem\mbf z_{k-1}\\ 0\eem+\wt\xi_k\wt{\mbf p}_k=\mbf x_k^{\rm L}+\wt\xi_k\wt{\mbf p}_k.$$ Recall that $\xi_2=\xi_4=\cdots=\xi_{2j}=0$ (see (\ref{xieven})).  It follows from (\ref{xi0}) that $\wt\xi_{2j}=0$. Hence, we have $$\mbf x_{2j}^{\rm G}=\mbf x_{2j}^{\rm L}.$$ We summarize the above theoretical results in the following theorem.

\begin{theorem}\label{theo11} Let $\mbf x_k^{\rm L}$ and $\mbf x_k^{\rm G}$ be the iterates generated at iteration $k$ of {\rm S$^3$LQ} and {\rm S$^3$CG}, respectively. As long as the algorithms have not yet converged, we have $\mbf x_{2j}^{\rm L}=\mbf x_{2j+1}^{\rm L}=\mbf x_{2j}^{\rm G}$ for $j\geq 1$.	
\end{theorem}

\subsection{Relation to USYMLQ and USYMQR} We will now show that S$^3$LQ and S$^3$MR are mathematically equivalent to two existing methods, USYMLQ and USYMQR \cite{saunders1988two}, for shifted skew-symmetric linear systems. 

Saunders--Simon--Yip tridiagonalization (Algorithm 3) applied to skew-symmetric $\mbf S$ and $\mbf c=\mbf b$ yields (see section 2.3) $$\mbf S\wt{\mbf V}_\ell=\wt{\mbf U}_\ell\wt{\mbf H}_\ell,\qquad \wt{\mbf U}_\ell=\wt{\mbf V}_\ell\wt{\mbf D}_\ell\mbf D_\ell,\qquad \wt{\mbf H}_\ell=\wt{\mbf U}_\ell^\top\mbf S\wt{\mbf V}_\ell=\bem 0 & \gamma_2 &&\\ \gamma_2 & 0 & \ddots &  \\ &\ddots &\ddots &\gamma_\ell  \\ &&\gamma_\ell & 0 \eem.$$ For $\mbf A=\alpha\mbf I+\mbf S$, we have $$\mbf A\wt{\mbf V}_\ell=\alpha\wt{\mbf V}_\ell+\mbf S\wt{\mbf V}_\ell=\alpha\wt{\mbf U}_\ell\mbf D_\ell\wt{\mbf D}_\ell+\wt{\mbf U}_\ell\wt{\mbf H}_\ell=\wt{\mbf U}_\ell\wt{\mbf T}_\ell,$$ where $$\wt{\mbf T}_\ell:=\alpha\mbf D_\ell\wt{\mbf D}_\ell+ \wt{\mbf H}_\ell=\bem \alpha & \gamma_2 &&\\ \gamma_2 & -\alpha & \ddots &  \\ &\ddots &\ddots & \gamma_\ell  \\ &&\gamma_\ell & (-1)^{\ell+1}\alpha \eem.$$

Let $\wt{\mbf T}_{k+1,k}$ denote the $(k+1)\times k$ upper-left submatrix of $\wt{\mbf T}_\ell$. The USYMQR points $\mbf x_k^{\rm USYMQR}$ satisfy $$\mbf x_k^{\rm USYMQR}=\wt{\mbf V}_k\wt{\mbf z}_k,\quad \wt{\mbf z}_k=\argmin_{\wt{\mbf z}\in\mbbr^k}\|\gamma_1\mbf e_1-\wt{\mbf T}_{k+1,k}\wt{\mbf z}\|_2.$$ Let $\wt{\mbf T}_{k+1,k}=\wt{\mbf Q}_{k+1,k}\wt{\mbf R}_k$ be the reduced QR factorization of $\wt{\mbf T}_{k+1,k}$. We have $\wt{\mbf z}_k=\wt{\mbf R}_k^{-1}\wt{\mbf Q}_{k+1,k}^\top\gamma_1\mbf e_1$. If $\mbf z_k=\wt{\mbf Q}_{k+1,k}^\top\gamma_1\mbf e_1$ and $\mbf P_k=\wt{\mbf V}_k\wt{\mbf R}_k^{-1}$, then $\mbf x_k^{\rm USYMQR}=\mbf P_k\mbf z_k$. From here on one can proceed as MINRES \cite{paige1975solut} for symmetric linear systems, and compute $\mbf z_k$ and the columns of $\mbf P_k$ by simple short recursions. For the case $\mbf S^\top=-\mbf S $ and $\mbf c=\mbf b$, we have (see section 2.3) $$\wt{\mbf V}_k=\mbf W_k\wt{\mbf D}_k,\quad \wt{\mbf H}_{k+1,k}=\mbf D_{k+1}\mbf H_{k+1,k}\wt{\mbf D}_k.$$ Then it holds that $\wt{\mbf T}_{k+1,k}=\mbf D_{k+1}\mbf T_{k+1,k}\wt{\mbf D}_k$. It is straightforward to show that $$\mbf x_k^{\rm USYMQR}=\wt{\mbf V}_k\wt{\mbf z}_k=\wt{\mbf V}_k\wt{\mbf D}_k\mbf y_k^{\rm M}=\mbf W_k\mbf y_k^{\rm M}=\mbf x_k^{\rm M}.$$ Therefore, S$^3$MR is mathematically equivalent to USYMQR for $\mbf A=\alpha\mbf I+\mbf S$ and $\mbf c=\mbf b$. Analogously, we can show that S$^3$LQ is mathematically equivalent to USYMLQ for the case $\mbf A=\alpha\mbf I+\mbf S$ and $\mbf c=\mbf b$.

\section{Numerical experiments}

The theoretical results of this work are obtained under the assumption that exact arithmetic is used. We illustrate our theoretical findings on artificially generated examples. Our goal here is to observe and confirm the theoretical properties (Propositions \ref{prop2}, \ref{prop3}, and \ref{leq}, Theorems \ref{theo7} and \ref{theo11}) of the various methods that we have studied.  

\begin{figure}[htb]
\centerline{\epsfig{figure=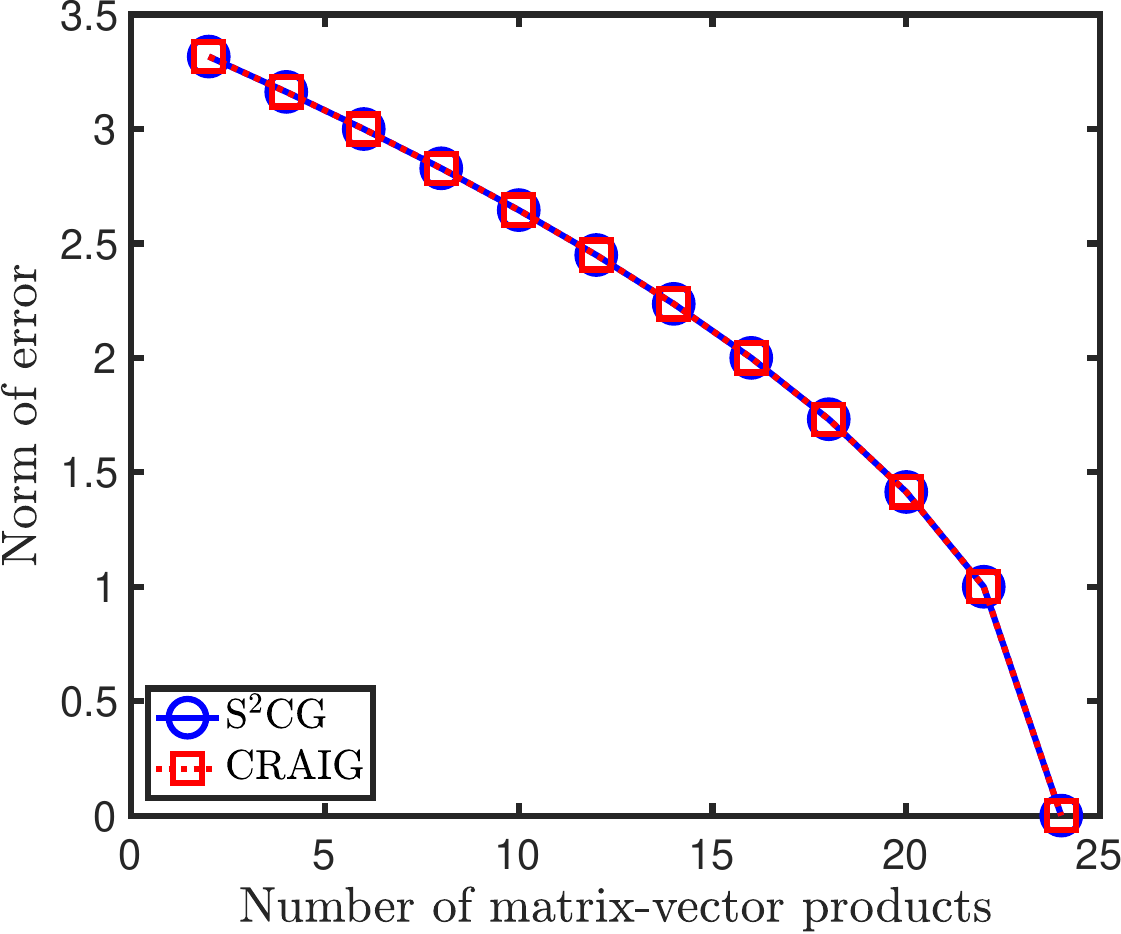,height=2.5in}\qquad\epsfig{figure=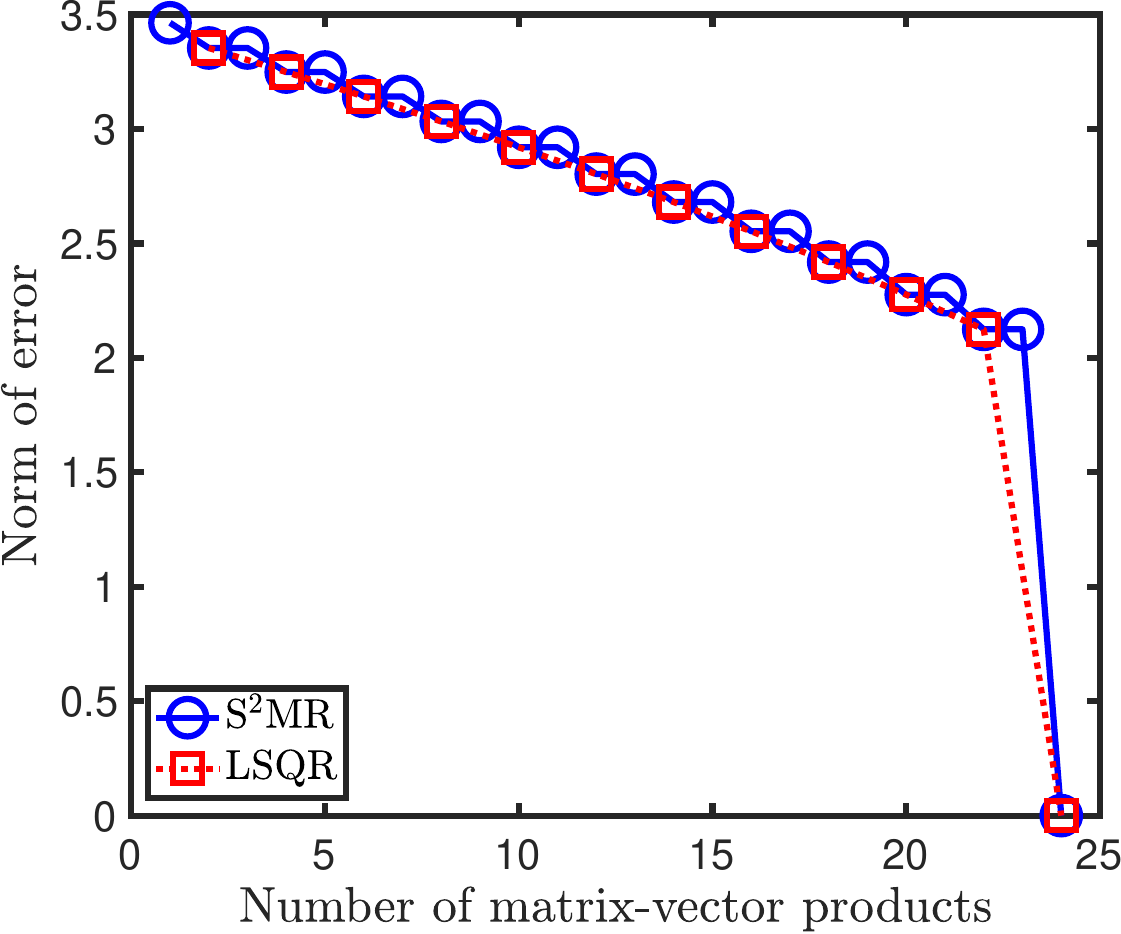,height=2.5in}}
\caption{Convergence history of S$^2$CG, CRAIG, S$^2$MR, and LSQR for a singular consistent skew-symmetric system. Left: S$^2$CG and CRAIG. Right: S$^2$MR and LSQR.}
\label{fig1}
\end{figure}

\subsection{Singular skew-symmetric linear systems}

We first consider the skew-symmetric consistent system $\mbf S\mbf x=\mbf b$ where \beq\label{ex1}\mbf S=\bem 0 & 1 &&\\ -1 & 0 & \ddots &  \\ &\ddots &\ddots & 1 \\ && -1 & 0 \eem\in\mbbr^{n\times n},\qquad \mbf b=\bem 1/\sqrt{2} \\ 0 \\ \vdots \\  0\\ -1/\sqrt{2}\eem\in\mbbr^n.\eeq This system is from Experiment 4.1 of \cite[section 4]{brown1997gmres}. We set $n=49$ and hence $\mbf S$ is singular. The pseudoinverse solution is $\mbf x_\star=\bem 0 &1/\sqrt{2}&0 & 1/\sqrt{2} &\cdots & 0 & 1/\sqrt{2}& 0 \eem^\top$. The errors of S$^2$CG, CRAIG, S$^2$MR, and LSQR are given in Figure \ref{fig1}. S$^2$CG and S$^2$MR safely terminate with $\mbf x_{24}^{\rm G}=\mbf x_\star$ and $\mbf x_{24}^{\rm M}=\mbf x_\star$. Our numerical results show that $\mbf x_{2j}^{\rm G}=\mbf x_j^{\rm CRAIG}$ and $\mbf x_{2j}^{\rm M}=\mbf x_{2j+1}^{\rm M}=\mbf x_j^{\rm LSQR}$ hold to within machine precision (see Propositions \ref{prop2} and \ref{prop3}). 
 
\begin{figure}[htb]
\centerline{\epsfig{figure=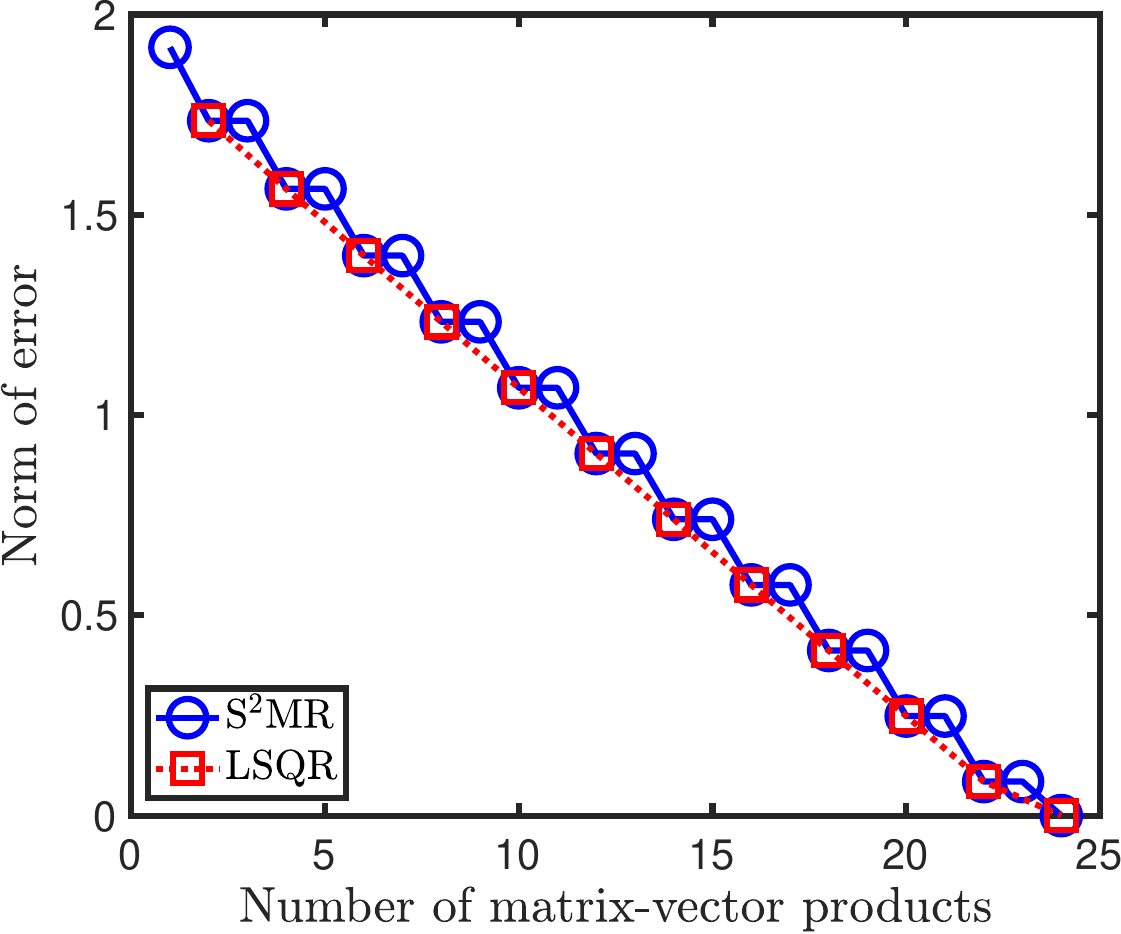,height=2.5in}}
\caption{Convergence history of S$^2$MR and LSQR for a singular inconsistent skew-symmetric system.}
\label{fig2}
\end{figure}

We then consider the singular inconsistent system $\mbf S\mbf x=\mbf b$ with the same $\mbf S$ as in (\ref{ex1}), but $\mbf b$ replaced with $\bem 1/\sqrt{2} & 0 & \cdots &  0& 1/\sqrt{2}\eem^\top$. The pseudoinverse solution is computed by using MATLAB's {\tt pinv} function. The errors of S$^2$MR and LSQR are given in Figure \ref{fig2}. S$^2$MR returns the pseudoinverse solution at the 24th step. We also observe that $\mbf x_{2j}^{\rm M}=\mbf x_{2j+1}^{\rm M}=\mbf x_j^{\rm LSQR}$ holds to within machine precision (see Proposition \ref{prop3}). 

\subsection{Shifted skew-symmetric linear systems}

\begin{figure}[htb]
\centerline{\epsfig{figure=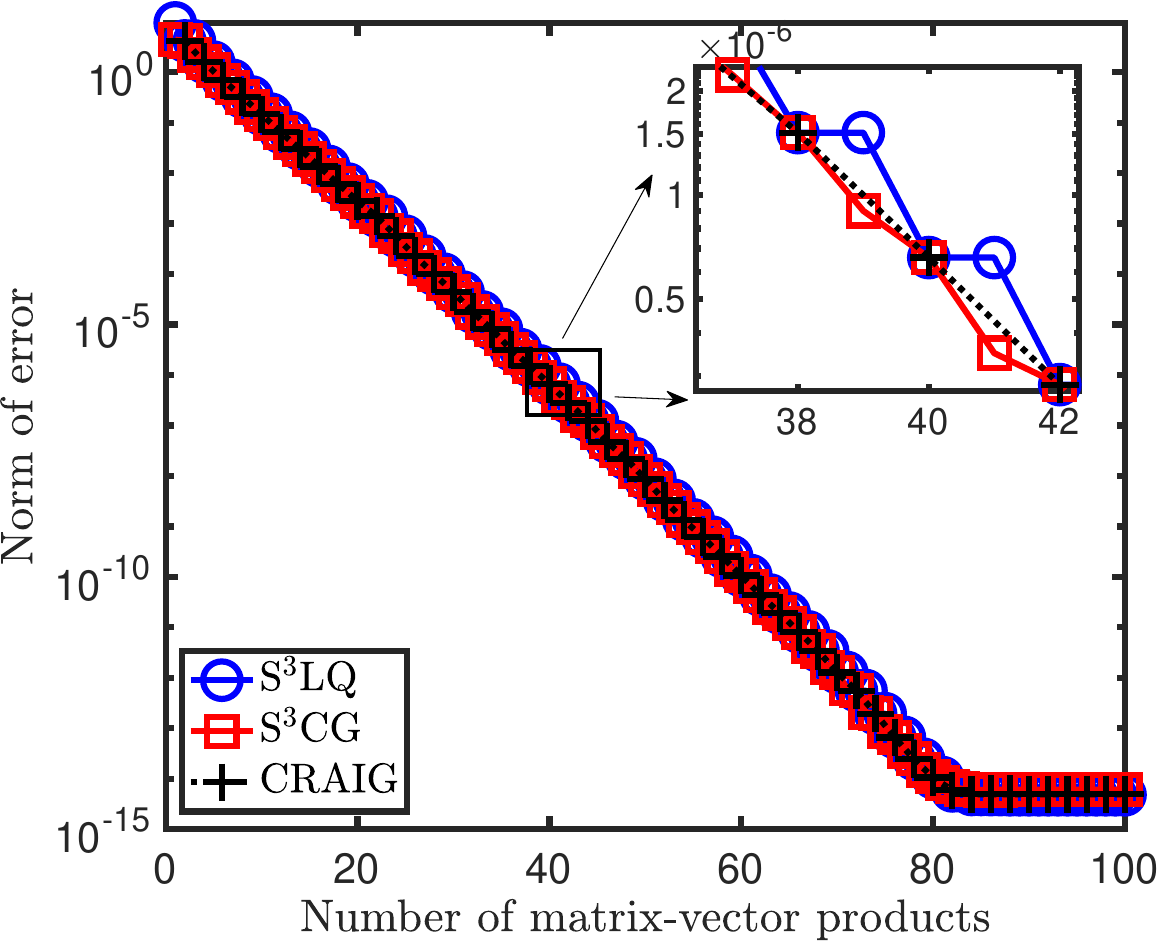,height=2.5in}\quad \epsfig{figure=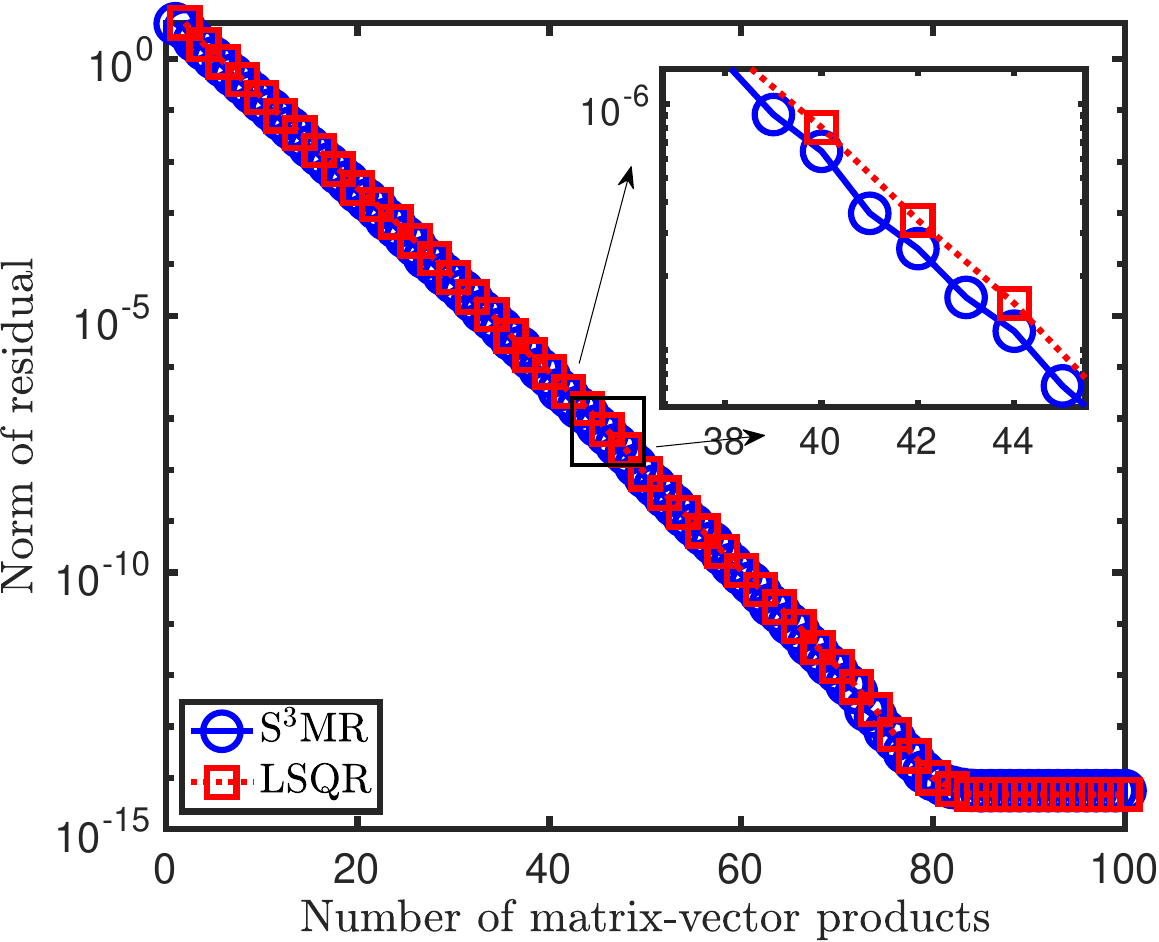,height=2.5in}}
\caption{Convergence history of S$^3$LQ, S$^3$CG, CRAIG, S$^3$MR, and LSQR for a shifted skew-symmetric system. Left: S$^3$LQ, S$^3$CG, and CRAIG. Right: S$^3$MR and LSQR.}
\label{fig3}
\end{figure}

The skew-symmetric matrix $\mbf S$ is obtained from a simple finite difference discretization of a constant convective term on the unit square in two dimensions, $(0,1)\times (0,1)$. For a given positive integer $m$ and a real scalar $\sigma$, we denote by $\mbf S_m(\sigma)$ the $m\times m$ skew-symmetric tridiagonal matrix $$\mbf S_m(\sigma)=\bem 0 & \sigma &&\\ -\sigma & 0 & \ddots &  \\ &\ddots &\ddots & \sigma \\ && -\sigma & 0 \eem\in\mbbr^{m\times m}.$$ The skew-symmetric matrix $\mbf S$ in $\mbf A=\alpha\mbf I+\mbf S$ is defined by $$\mbf S=\mbf I_m\otimes \mbf S_m(\sigma_1)+\mbf S_m(\sigma_2)\otimes\mbf I_m,$$ where $\otimes$ denotes the Kronecker product, and $\sigma_1$ and $\sigma_2$ are preselected scalars. Evidently $\mbf S$ and $\mbf A$ are of dimensions $m^2\times m^2$. The right-hand side is randomly generated. We show results with $m=15$, $\alpha=0.8$, $\sigma_1=0.4$ and $\sigma_2=0.6$. The errors of S$^3$LQ, S$^3$CG, and CRAIG, and the residuals of S$^3$MR and LSQR are given in Figure \ref{fig3}. Our numerical results show that $\mbf x_{2j}^{\rm L}=\mbf x_{2j+1}^{\rm L}=\mbf x_{2j}^{\rm G}=\mbf x_j^{\rm CRAIG}$ hold to within machine precision (see Theorems \ref{theo7} and \ref{theo11}). It is also observed that $\|\mbf b-\mbf A\mbf x_{2j}^{\rm M}\|_2<\|\mbf b-\mbf A\mbf x_j^{\rm LSQR}\|_2$ (see Proposition \ref{leq}).

\section{Concluding remarks}

We have addressed the performance of S$^2$CG and S$^2$MR on skew-symmetric linear systems $\mbf A\mbf x=\mbf b$ when $\mbf A$ is singular. The theoretical results indicate that S$^2$CG can find the pseudoinverse solution of consistent systems, while S$^2$MR can find the pseudoinverse solution of consistent and inconsistent systems.

S$^3$LQ is an iterative method for shifted skew-symmetric linear systems, with the property that it ensures monotonic reduction in the error $\|\mbf x_k^{\rm L}-\mbf A^{-1}\mbf b\|_2$. In deriving it we have completed the triad of solvers S$^3$CG, S$^3$MR, S$^3$LQ for solving shifted skew-symmetric linear systems based on Lanczos tridiagonalization. They correspond to CG, MINRES, and SYMMLQ for solving symmetric linear systems, respectively. Our implementation of S$^3$CG has the advantage that it emphasizes how S$^3$CG is related to CG. We have given the relation among iterations of S$^3$LQ, S$^3$CG, and CRAIG, namely, $\mbf x_{2j}^{\rm L}=\mbf x_{2j+1}^{\rm L}=\mbf x_{2j}^{\rm G}=\mbf x_j^{\rm CRAIG}$, and the relation between residuals of S$^3$MR and LSQR, namely, $\|\mbf b-\mbf A\mbf x_{2k}^{\rm M}\|_2\leq\|\mbf b-\mbf A\mbf x_k^{\rm LSQR}\|_2$. By the relation between Lanczos tridiagonalization and Saunders--Simon--Yip tridiagonalization, we have shown that S$^3$LQ and S$^3$MR are mathematically equivalent to USYMLQ and USYMQR for shifted skew-symmetric linear systems, respectively.

Practical computations always incorporate preconditioning. We mention that some preconditioning techniques have been briefly discussed by Greif and Varah \cite{greif2009itera}, Greif et al. \cite{greif2016numer}, and Manguo\u{g}lu and Mehrmann \cite{manguoglu2021two}. We also notice that a structured shifted skew-symmetric problem has been considered by Greif \cite{greif2022struc}. We leave preconditioning techniques and structured shifted problems for future research. Note also that in this paper we assume only exact arithmetic. Effects of finite precision arithmetic will be discussed in the future.

\appendix
\section{Golub--Kahan bidiagonalization for shifted skew-symmetric matrices} Golub--Kahan bidiagonalization (Algorithm 2) for shifted skew-symmetric $\mbf A=\alpha\mbf I+\mbf S$ ($\alpha\neq 0$) and nonzero $\mbf b$ must stop in $\ell_0=\lceil\ell/2\rceil$ steps with $\beta_{\ell_0+1}=0$ because $\mbf b\in\ran(\mbf A)$ always holds. We have $$\mbf A\mbf V_{\ell_0}=\mbf U_{\ell_0}\mbf B_{\ell_0},\qquad \mbf A^\top\mbf U_{\ell_0}=\mbf V_{\ell_0}\mbf B_{\ell_0}^\top,\qquad \beta_1=\gamma_1,\qquad \alpha_1=\sqrt{\alpha^2+\gamma_2^2}.$$

\begin{proposition} For $j\geq 1$, we have $\alpha_j>\gamma_{2j}$ and $\beta_{j+1}= \gamma_{2j+1}\gamma_{2j}/\alpha_j<\gamma_{2j+1}$. \end{proposition} 
\proof Golub--Kahan bidiagonalization (Algorithm 2) for $\mbf S$ and $\mbf b$ gives (see section 2.2) \begin{align*} \gamma_{2j+1}(-1)^j\mbf w_{2j+1} &=\mbf S(-1)^j\mbf w_{2j} - \gamma_{2j}(-1)^{j-1}\mbf w_{2j-1}, \\ \gamma_{2j+2}(-1)^{j+1}\mbf w_{2j+2} &=-\mbf S(-1)^j\mbf w_{2j+1}-\gamma_{2j+1}(-1)^j\mbf w_{2j}.\end{align*} Golub--Kahan bidiagonalization (Algorithm 2) for $\mbf A=\alpha\mbf I+\mbf S$ and $\mbf b$ gives \begin{align*} \beta_{j+1}\mbf u_{j+1} &=(\alpha\mbf I+\mbf S)\mbf v_j - \alpha_j\mbf u_j, \\ \alpha_{j+1}\mbf v_{j+1} &=(\alpha\mbf I-\mbf S)\mbf u_{j+1}-\beta_{j+1}\mbf v_j.\end{align*} By induction on $j$, it is straightforward to show that $$\mbf u_{j+1}=(-1)^j\mbf w_{2j+1},\qquad \mbf v_{j+1}\in\spa\{\mbf w_1,\cdots,\mbf w_{2j+2}\}.$$ Direct calculations yield $$\alpha_{j+1}\mbf w_{2j+1}^\top\mbf v_{j+1}=\mbf w_{2j+1}^\top(\alpha\mbf I-\mbf S)\mbf u_{j+1}-\beta_{j+1}\mbf w_{2j+1}^\top\mbf v_j=(-1)^j\alpha$$ and \begin{align*}\alpha_{j+1}\mbf w_{2j+2}^\top\mbf v_{j+1}&=\mbf w_{2j+2}^\top(\alpha\mbf I-\mbf S)\mbf u_{j+1}-\beta_{j+1}\mbf w_{2j+2}^\top\mbf v_j=-\mbf w_{2j+2}^\top\mbf S(-1)^j\mbf w_{2j+1}\\ & = \mbf w_{2j+2}^\top(\gamma_{2j+2}(-1)^{j+1}\mbf w_{2j+2}+\gamma_{2j+1}(-1)^j\mbf w_{2j})\\ &=(-1)^{j+1}\gamma_{2j+2}.\end{align*} Therefore, we have $$\|\mbf v_{j+1}\|_2^2=1\geq \frac{(-1)^{2j}\alpha^2+(-1)^{2j+2}\gamma_{2j+2}^2}{\alpha_{j+1}^2},$$ which gives $$\alpha_{j+1}\geq\sqrt{\alpha^2+\gamma_{2j+2}^2}>\gamma_{2j+2}.$$ Also by direct calculations, we obtain \begin{align*}\beta_{j+1}=\mbf u_{j+1}^\top\beta_{j+1}\mbf u_{j+1} &=(-1)^j\mbf w_{2j+1}^\top(\alpha\mbf I+\mbf S)\mbf v_j-\mbf u_{j+1}^\top\alpha_j\mbf u_j\\ &=(-1)^j\mbf w_{2j+1}^\top\mbf S\mbf v_j=(-1)^j(-\mbf S\mbf w_{2j+1})^\top\mbf v_j\\ &=(-1)^j(\gamma_{2j+1}\mbf w_{2j}-\gamma_{2j+2}\mbf w_{2j+2})^\top\mbf v_j\\ &=(-1)^j\gamma_{2j+1}\mbf w_{2j}^\top\mbf v_j\\ & = \gamma_{2j+1}\gamma_{2j}/\alpha_j< \gamma_{2j+1}.\end{align*} This completes the proof. \endproof

\section*{Funding/Conflicts of interests/Competing interests} This work was supported by the National Natural Science Foundation of China (No.12171403 and No.11771364), the Natural Science Foundation of Fujian Province of China (No.2020J01030), and the Fundamental Research Funds for the Central Universities (No.20720210032). There are no conflicts of interests to this work. The authors have no competing interests to declare that are relevant to the content of this article.		
	
{\small 

\begin{thebibliography}{10}

\bibitem{arnoldi1951princ}
W.~E. Arnoldi.
\newblock The principle of minimized iteration in the solution of the matrix
  eigenvalue problem.
\newblock {\em Quart. Appl. Math.}, 9:17--29, 1951.

\bibitem{bai2003hermi}
Z.-Z. Bai, G.~H. Golub, and M.~K. Ng.
\newblock Hermitian and skew-{H}ermitian splitting methods for non-{H}ermitian
  positive definite linear systems.
\newblock {\em SIAM J. Matrix Anal. Appl.}, 24(3):603--626, 2003.

\bibitem{brown1997gmres}
P.~N. Brown and H.~F. Walker.
\newblock G{MRES} on (nearly) singular systems.
\newblock {\em SIAM J. Matrix Anal. Appl.}, 18(1):37--51, 1997.

\bibitem{calvetti2000gmres}
D.~Calvetti, B.~Lewis, and L.~Reichel.
\newblock G{MRES}-type methods for inconsistent systems.
\newblock {\em Linear Algebra Appl.}, 316(1-3):157--169, 2000.

\bibitem{cao2002note}
Z.-H. Cao and M.~Wang.
\newblock A note on {K}rylov subspace methods for singular systems.
\newblock {\em Linear Algebra Appl.}, 350:285--288, 2002.

\bibitem{chen2007regul}
J.~Chen and Z.~Shen.
\newblock Regularized conjugate gradient method for skew-symmetric and
  indefinite system of linear equations and applications.
\newblock {\em Appl. Math. Comput.}, 187(2):1484--1494, 2007.

\bibitem{chen2023accel}
L.~Chen and J.~Wei.
\newblock Accelerated gradient and skew-symmetric splitting methods for a class
  of monotone operator equations.
\newblock arXiv preprint arXiv:2303.09009, 2023.

\bibitem{choi2013minim}
S.-C.~T. Choi.
\newblock Minimal residual methods for complex symmetric, skew symmetric, and
  skew {H}ermitian systems.
\newblock arXiv preprint arXiv:1304.6782, 2013.

\bibitem{choi2011minre}
S.-C.~T. Choi, C.~C. Paige, and M.~A. Saunders.
\newblock M{INRES}-{QLP}: {A} {K}rylov subspace method for indefinite or
  singular symmetric systems.
\newblock {\em SIAM J. Sci. Comput.}, 33(4):1810--1836, 2011.

\bibitem{concus1976gener}
P.~Concus and G.~H. Golub.
\newblock A generalized conjugate gradient method for nonsymmetric systems of
  linear equations.
\newblock In {\em Computing methods in applied sciences and engineering
  ({S}econd {I}nternat. {S}ympos., {V}ersailles, 1975), {P}art 1}, Lecture
  Notes in Econom. and Math. Systems, Vol. 134, pages 56--65. Springer, Berlin,
  1976.

\bibitem{craig1955nstep}
E.~J. Craig.
\newblock The {$N$}-step iteration procedures.
\newblock {\em J. Math. and Phys.}, 34:64--73, 1955.

\bibitem{demmel1997appli}
J.~W. Demmel.
\newblock {\em Applied Numerical Linear Algebra}.
\newblock Society for Industrial and Applied Mathematics (SIAM), Philadelphia,
  PA, 1997.

\bibitem{du2017any}
K.~Du, J.~Duintjer~Tebbens, and G.~Meurant.
\newblock Any admissible harmonic {R}itz value set is possible for {GMRES}.
\newblock {\em Electron. Trans. Numer. Anal.}, 47:37--56, 2017.

\bibitem{eisenstat1983note}
S.~C. Eisenstat.
\newblock A note on the generalized conjugate gradient method.
\newblock {\em SIAM J. Numer. Anal.}, 20(2):358--361, 1983.

\bibitem{eisenstat2015equiv}
S.~C. Eisenstat.
\newblock Equivalence of {K}rylov subspace methods for skew-symmetric linear
  systems.
\newblock arXiv preprint arXiv:1512.0031, 2015.

\bibitem{estrin2019lnlq}
R.~Estrin, D.~Orban, and M.~A. Saunders.
\newblock L{NLQ}: {A}n iterative method for least-norm problems with an error
  minimization property.
\newblock {\em SIAM J. Matrix Anal. Appl.}, 40(3):1102--1124, 2019.

\bibitem{estrin2019lslq}
R.~Estrin, D.~Orban, and M.~A. Saunders.
\newblock L{SLQ}: {A}n iterative method for linear least-squares with an error
  minimization property.
\newblock {\em SIAM J. Matrix Anal. Appl.}, 40(1):254--275, 2019.

\bibitem{fan2024itera}
J.-J. Fan.
\newblock On iterative methods for skew-symmetric and shifted skew-symmetric
  linear systems.
\newblock Master's thesis, Xiamen University, 2024 (in preparation).

\bibitem{fong2011lsmr}
D.~C.-L. Fong and M.~Saunders.
\newblock L{SMR}: {A}n iterative algorithm for sparse least-squares problems.
\newblock {\em SIAM J. Sci. Comput.}, 33(5):2950--2971, 2011.

\bibitem{golub1965calcu}
G.~Golub and W.~Kahan.
\newblock Calculating the singular values and pseudo-inverse of a matrix.
\newblock {\em SIAM J. Numer. Anal.}, 2(2):205--224, 1965.

\bibitem{greenbaum1997itera}
A.~Greenbaum.
\newblock {\em Iterative {M}ethods for {S}olving {L}inear {S}ystems}, volume~17
  of {\em Frontiers in Applied Mathematics}.
\newblock Society for Industrial and Applied Mathematics (SIAM), Philadelphia,
  PA, 1997.

\bibitem{greif2022struc}
C.~Greif.
\newblock Structured shifts for skew-symmetric matrices.
\newblock {\em Electron. Trans. Numer. Anal.}, 55:455--468, 2022.

\bibitem{greif2016numer}
C.~Greif, C.~C. Paige, D.~Titley-Peloquin, and J.~M. Varah.
\newblock Numerical equivalences among {K}rylov subspace algorithms for
  skew-symmetric matrices.
\newblock {\em SIAM J. Matrix Anal. Appl.}, 37(3):1071--1087, 2016.

\bibitem{greif2009itera}
C.~Greif and J.~M. Varah.
\newblock Iterative solution of skew-symmetric linear systems.
\newblock {\em SIAM J. Matrix Anal. Appl.}, 31(2):584--601, 2009.

\bibitem{gu2009skew}
C.~Gu and H.~Qian.
\newblock Skew-symmetric methods for nonsymmetric linear systems with multiple
  right-hand sides.
\newblock {\em J. Comput. Appl. Math.}, 223(2):567--577, 2009.

\bibitem{guducu2022nonhe}
C.~G\"{u}d\"{u}c\"{u}, J.~Liesen, V.~Mehrmann, and D.~B. Szyld.
\newblock On non-{H}ermitian positive (semi)definite linear algebraic systems
  arising from dissipative {H}amiltonian {DAE}s.
\newblock {\em SIAM J. Sci. Comput.}, 44(4):A2871--A2894, 2022.

\bibitem{hageman1980equiv}
L.~A. Hageman, F.~T. Luk, and D.~M. Young.
\newblock On the equivalence of certain iterative acceleration methods.
\newblock {\em SIAM J. Numer. Anal.}, 17(6):852--873, 1980.

\bibitem{hayami2011geome}
K.~Hayami and M.~Sugihara.
\newblock A geometric view of {K}rylov subspace methods on singular systems.
\newblock {\em Numer. Linear Algebra Appl.}, 18(3):449--469, 2011.

\bibitem{hestenes1952metho}
M.~R. Hestenes and E.~Stiefel.
\newblock Methods of conjugate gradients for solving linear systems.
\newblock {\em J. Research Nat. Bur. Standards}, 49:409--436 (1953), 1952.

\bibitem{huang1999itera}
Y.~Huang, A.~J. Wathen, and L.~Li.
\newblock An iterative method for skew-symmetric systems.
\newblock {\em Information}, 2(2):147--153, 1999.

\bibitem{idema2007minim}
R.~Idema and C.~Vuik.
\newblock A minimal residual method for shifted skew-symmetric systems.
\newblock Technical Report 07-09, Delft University of Technology, Delft, The
  Netherlands, 2007.

\bibitem{idema2023compa}
R.~Idema and C.~Vuik.
\newblock A comparison of {K}rylov methods for shifted skew-symmetric systems.
\newblock arXiv preprint arXiv:2304.04092, 2023.

\bibitem{ipsen1998idea}
I.~C.~F. Ipsen and C.~D. Meyer.
\newblock The idea behind {K}rylov methods.
\newblock {\em Amer. Math. Monthly}, 105(10):889--899, 1998.

\bibitem{jiang2007algor}
E.~Jiang.
\newblock Algorithm for solving shifted skew-symmetric linear system.
\newblock {\em Front. Math. China}, 2(2):227--242, 2007.

\bibitem{kammerer1972conve}
W.~J. Kammerer and M.~Z. Nashed.
\newblock On the convergence of the conjugate gradient method for singular
  linear operator equations.
\newblock {\em SIAM J. Numer. Anal.}, 9:165--181, 1972.

\bibitem{lanczos1950itera}
C.~Lanczos.
\newblock An iteration method for the solution of the eigenvalue problem of
  linear differential and integral operators.
\newblock {\em J. Research Nat. Bur. Standards}, 45:255--282, 1950.

\bibitem{manguoglu2021two}
M.~Manguo\u{g}lu and V.~Mehrmann.
\newblock A two-level iterative scheme for general sparse linear systems based
  on approximate skew-symmetrizers.
\newblock {\em Electron. Trans. Numer. Anal.}, 54:370--391, 2021.

\bibitem{meurant2014neces}
G.~Meurant.
\newblock Necessary and sufficient conditions for {GMRES} complete and partial
  stagnation.
\newblock {\em Appl. Numer. Math.}, 75:100--107, 2014.

\bibitem{morikuni2018gmres}
K.~Morikuni and M.~Rozlo\v{z}n\'{\i}k.
\newblock On {GMRES} for singular {EP} and {GP} systems.
\newblock {\em SIAM J. Matrix Anal. Appl.}, 39(2):1033--1048, 2018.

\bibitem{paige1974bidia}
C.~C. Paige.
\newblock Bidiagonalization of matrices and solutions of the linear equations.
\newblock {\em SIAM J. Numer. Anal.}, 11:197--209, 1974.

\bibitem{paige1975solut}
C.~C. Paige and M.~A. Saunders.
\newblock Solutions of sparse indefinite systems of linear equations.
\newblock {\em SIAM J. Numer. Anal.}, 12(4):617--629, 1975.

\bibitem{paige1982lsqr}
C.~C. Paige and M.~A. Saunders.
\newblock L{SQR}: {A}n algorithm for sparse linear equations and sparse least
  squares.
\newblock {\em ACM Trans. Math. Software}, 8(1):43--71, 1982.

\bibitem{rapoport1978nonli}
D.~Rapoport.
\newblock {\em A {N}onlinear {L}anczos {A}lgorithm and the {S}tationary
  {N}avier--{S}tokes {E}quation}.
\newblock ProQuest LLC, Ann Arbor, MI, 1978.
\newblock Thesis (Ph.D.)--New York University.

\bibitem{reichel2005break}
L.~Reichel and Q.~Ye.
\newblock Breakdown-free {GMRES} for singular systems.
\newblock {\em SIAM J. Matrix Anal. Appl.}, 26(4):1001--1021, 2005.

\bibitem{saad2003itera}
Y.~Saad.
\newblock {\em Iterative Methods for Sparse Linear Systems}.
\newblock Society for Industrial and Applied Mathematics, Philadelphia, PA,
  second edition, 2003.

\bibitem{saad1986gmres}
Y.~Saad and M.~H. Schultz.
\newblock G{MRES}: a generalized minimal residual algorithm for solving
  nonsymmetric linear systems.
\newblock {\em SIAM J. Sci. Statist. Comput.}, 7(3):856--869, 1986.

\bibitem{saunders1988two}
M.~A. Saunders, H.~D. Simon, and E.~L. Yip.
\newblock Two conjugate-gradient-type methods for unsymmetric linear equations.
\newblock {\em SIAM J. Numer. Anal.}, 25(4):927--940, 1988.

\bibitem{sidi1999unifi}
A.~Sidi.
\newblock A unified approach to {K}rylov subspace methods for the
  {D}razin-inverse solution of singular nonsymmetric linear systems.
\newblock {\em Linear Algebra Appl.}, 298(1-3):99--113, 1999.

\bibitem{smoch1999some}
L.~Smoch.
\newblock Some results about {GMRES} in the singular case.
\newblock {\em Numer. Algorithms}, 22(2):193--212, 1999.

\bibitem{smoch2007spect}
L.~Smoch.
\newblock Spectral behaviour of {GMRES} applied to singular systems.
\newblock {\em Adv. Comput. Math.}, 27(2):151--166, 2007.

\bibitem{szyld1993varia}
D.~B. Szyld and O.~B. Widlund.
\newblock Variational analysis of some conjugate gradient methods.
\newblock {\em East-West J. Numer. Math.}, 1:51--74, 1993.

\bibitem{trefethen1997numer}
L.~N. Trefethen and D.~Bau, III.
\newblock {\em Numerical Linear Algebra}.
\newblock Society for Industrial and Applied Mathematics (SIAM), Philadelphia,
  PA, 1997.

\bibitem{wei2000conve}
Y.~Wei and H.~Wu.
\newblock Convergence properties of {K}rylov subspace methods for singular
  linear systems with arbitrary index.
\newblock {\em J. Comput. Appl. Math.}, 114(2):305--318, 2000.

\bibitem{widlund1978lancz}
O.~Widlund.
\newblock A {L}anczos method for a class of nonsymmetric systems of linear
  equations.
\newblock {\em SIAM J. Numer. Anal.}, 15(4):801--812, 1978.

\end{thebibliography}

}
\end{document}